\documentclass[12pt]{article}
\usepackage{amsfonts}
\usepackage{amssymb}
\usepackage{polski}

\input{tcilatex}
\begin{document}

\begin{center}
\textbf{General Regular Variation, Popa Groups and Quantifier Weakening}

\bigskip \textbf{by \\[0pt]
N. H. Bingham and A. J. Ostaszewski}\\[0pt]

\bigskip

\textit{To Harry I.\ Miller (22 Feb 1939 - 15/16 December 2018), man and
mathematician, who died with his boots on}

\bigskip

`The soldiers' music and the rites of war

Speak loudly for him.' (Shakespeare, Hamlet Act V.2)

\bigskip
\end{center}

\noindent \textbf{Abstract.} We introduce \textit{general regular variation}%
, a theory of regular variation containing the existing Karamata,
Bojanic-Karamata/de Haan and Beurling theories as special cases. The
unifying theme is the Popa groups of our title viewed as locally compact
abelian ordered topological groups, together with their Haar measure and
Fourier theory. The power of this unified approach is shown by the
simplification it brings to the whole area of quantifier weakening, so
important in this field.

\bigskip


\bigskip

\noindent \textbf{Classification}: 26A03, 26A12, 33B99, 39B22.

\bigskip

\section{Introduction}

We recall the definition of \textit{Beurling slowly varying }functions $%
\varphi $ (see e.g. [BinGT \S\ 2.11], [BinO7]): these are positive,
measurable or Baire (i.e. have the Baire property, BP), are $o(x)$ at
infinity (or $O(x)$, depending on context), and, with%
\[
x\circ _{\varphi }t:=x+t\varphi (x)
\]%
the \textit{Popa }(or \textit{circle}) \textit{operation} (\S\ 2 below),
satisfy\textit{\ }%
\begin{equation}
\log \varphi (x\circ _{\varphi }t)-\log \varphi (x)\rightarrow 0:\qquad
\varphi (x\circ _{\varphi }t)/\varphi (x)\rightarrow 1.  \tag{$B$}
\end{equation}%
Such $\varphi $ will play the role of auxiliary functions below. For a
suitable auxiliary function $h$ and limit function $K,$ called the \textit{%
kernel, }consider also the limit relationship
\begin{equation}
\lbrack f(x\circ _{\varphi }t)-f(x)]/h(x)\rightarrow K(t),  \tag{$G$}
\end{equation}%
where $f$ here is the function of primary interest (`G for Goldie, G for
general': see e.g. [BinO6,10,11], [Ost2]). Specialising to $\varphi \equiv
1,h\equiv 1$ gives%
\begin{equation}
f(x+t)-f(x)\rightarrow K(t)  \tag{$K$}
\end{equation}%
(`K for Karamata'). This is the defining relationship for \textit{Karamata
regular variation }written additively (see e.g. [BinGT Ch. 1-3]: one needs
to be able to pass between the additive notation above, and the original
multiplicative notation, using the familiar exp-log isomorphism between the
additive group of reals (Haar measure Lebesgue measure) and the
multiplicative group of positive reals (Haar measure $\mathrm{d}x/x$).
Specialising instead to $\varphi \equiv 1,h$ slowly varying (in Karamata's
sense: [BinGT, Ch. 1]) gives

\begin{equation}
\lbrack f(x+t)-f(x)]/h(x)\rightarrow K(t),  \tag{$BKdH$}
\end{equation}%
the defining relationship for \textit{Bojanic-Karamata/de Haan regular
variation }[BinGT, Ch. 3], while specialising to $f$ = $\log \varphi ,$ $h$ $%
=1$, $K=0$ gives Beurling slow variation as above. We call the limit
relationship ($G$) above \textit{general regular variation, }as it contains
the other three. Below we give a unified treatment, using the
algebraicization provided by the \textit{Popa groups }of \S 2 below.

As usual (see e.g. [BinO1,9]), we pass between the measurable and Baire
cases (in any form of regular variation) `bitopologically' -- by passing
between the Euclidean and density topologies. The same will be true in the
Popa groups below, which are isomorphic to the reals algebraically and
bitopologically; we thus extend the terms Euclidean and density topologies
to these Popa isomorphs also.

\section{Popa groups}

Above we have used the Popa operation as a simplifying notational device for
the regular variation above (general or otherwise), involving limits as $%
x\rightarrow \infty $. But its usefulness is far greater, and is not
confined to limits, as emerged in [BinO7], [Ost1]. Here one allows other
auxiliary functions $h,$ with corresponding circle operations $\circ _{h}$.
This is most useful when the circle operation is associative, and this
requires $h$ to satisfy the\textit{\ Go\l \k{a}b-Schintzel functional
equation}:%
\begin{equation}
h(s\circ _{h}t)=h(s+h(s)t)=h(s)h(t)  \tag{$GS$}
\end{equation}%
(cf. [Jav]). Thus ($GS$) expresses homomorphy in this context, which will
occur in the regular variation context after the passage to the limit $%
x\rightarrow \infty $. Indeed, such an $h$ generates group structures on
subsets of $\mathbb{R}$, that are in fact isomorphic to the group $(\mathbb{R%
}_{+},\times ).$ It is to these \textit{Popa groups }[Pop] that we now turn.

Write $GS$ for the set of \textit{positive} solutions $h$ of ($GS$). It
emerges that, being thus bounded below, they are continuous and of the form%
\[
\eta (t)=\eta _{\rho }(t):=\eta (t)=1+\rho t
\]%
for $t>-1/\rho ,$ with the parameter $\rho \geq 0;$ for a proof see [Brz2]
and [BrzM], or the more direct [Ost3, \S 5] -- see also [AczD] and the
surveys [Brz1] and [Jab5]; cf. [Jab2], [Ost1]. For $\eta \in GS,$ put%
\[
\mathbb{G}_{\eta }^{\ast }:=\{x\in \mathbb{R}:\eta (x)\neq 0\}.
\]

Equipped with $\circ _{\eta },$ this is a group . When $\eta =\eta _{\rho }$
this operation is given explicitly by%
\[
x\circ _{\rho }y=x+y(1+\rho x),
\]%
so that $\mathbb{G}_{\rho }^{\ast }=\{x\in \mathbb{R}:x\neq \rho ^{\ast }\},$
where $\rho ^{\ast }=-1/\rho ,$ the \textit{Popa centre}. We interpret this
to mean $\rho ^{\ast }=-\infty $ for $\rho =0$ and to mean $\rho ^{\ast }=0$
for $\rho =+\infty .$

The operation $\circ _{\rho }$ may also be rendered by reference to the
equation ($GS$) in the current context:%
\[
\eta _{\rho }(x\circ _{\rho }y)=\eta _{\rho }(x)\eta _{\rho }(y)\qquad
(x,y\in \mathbb{G}_{\rho }),
\]%
and thereby to the underlying role of the multiplicative positive reals $%
\mathbb{R}_{+}$:
\[
x\circ _{\rho }y=\eta _{\rho }^{-1}(\eta _{\rho }(x)\eta _{\rho
}(x))=[(1+\rho x)(1+\rho y)-1]/\rho
\]%
(which gives for $\rho =1$ the \textit{circle operation }of ring theory: cf.
[Ost2, \S 2.1]). It emerges from here that $($except for the case $\rho =0$
where $\rho ^{\ast }=-\infty $ so that $\mathbb{G}_{0}^{\ast }=\mathbb{R}$)
the following subgroups of $\mathbb{G}_{\rho }^{\ast }$ are of greater
significance:%
\[
\mathbb{G}_{\rho }:=\{x\in \mathbb{R}:1+\rho x>0\}:\quad \quad \mathbb{G}%
_{\eta }:=\{x\in \mathbb{R}:\eta (x)>0\},
\]%
by virtue of being isomorphic with $(\mathbb{R}_{+},\times )$ when $\rho >0.$
(Likewise, the groups $\mathbb{G}_{\rho }^{\ast }$ are all isomorphic with $(%
\mathbb{R}^{\ast },\times ),$ with $\mathbb{R}^{\ast }:=\mathbb{R}\backslash
\{0\}.$)

As $\rho ^{\ast }=0$ for $\rho =+\infty ,$ the group $(\mathbb{R}_{+},\times
)$ may itself be viewed conveniently as $\mathbb{G}_{\infty },$ or perhaps
more accurately as the rescaled limit of $\mathbb{G}_{\rho }$ as $\rho
\rightarrow \infty ,$ as follows:
\[
(x\circ _{\rho }y)/\rho =[(1+\rho x)(1+\rho x)-1]/\rho ^{2}=x/\rho +y/\rho
+xy\rightarrow xy,\text{ as }\rho \rightarrow \infty .
\]

We note that one has $1_{\mathbb{G}}=0$ for $\mathbb{G=G}_{\rho }\mathbb{\ }$%
except for $\mathbb{G=G}_{\infty },$ when $1_{\mathbb{G}}=1.$ The \textit{%
inverse} of $t$ in $\mathbb{G}_{\eta }$ will be denoted by $t_{\eta }^{-1}$
(or $t_{\rho }^{-1},$ if more convenient); here

\[
t_{\eta }^{-1}=-t/\eta (t).
\]

We will also need to designate location to either side of $1_{\mathbb{G}}=0$%
, using the notation%
\[
\mathbb{G}_{\rho }^{+}:=\{x\in \mathbb{G}_{\rho }:x>0\text{ }\&\text{ }%
1+\rho x>0\}
\]%
and $\mathbb{G}_{\rho }^{-}:=\{x\in \mathbb{G}_{\rho }:x<0$ $\&$ $1+\rho
x>0\}.$

Viewing the Popa operation as a conjugacy via the isomorphism $\eta _{\rho
}, $%
\begin{equation}
x\circ _{\rho }y=[(1+\rho x)(1+\rho y)-1]/\rho ,  \tag{conj}
\end{equation}%
demonstrates that $\circ _{\rho }$ may be expressed in terms of the \textit{%
ring operations }of $\mathbb{R}$, and so permits other features of $\mathbb{R%
}$ to be imported into $\mathbb{G}_{\rho }$. There are several possibilities
here. The Popa groups may inherit either of the two canonical topological
structures of their isomorphs, again enabling bitopological passage between
them (as in \S 1). Thus they inherit a \textit{Euclidean} \textit{topology},
from which they derive their own metric structures; this is generated by
(open) intervals, and makes $\mathbb{G}_{\rho }$ a \textit{locally compact}
abelian topological group. In turn this allows reference to \textit{Haar
measure}, and so to the second possibility: the \textit{Haar-density}
topology of $\mathbb{G}_{\rho },$ which which agrees with the topology
induced on $\mathbb{G}_{\rho }$ by the (Lebesgue) \textit{density} \textit{%
topology} on $\mathbb{R}$ (corresponding to \textit{Lebesgue measure} $%
\lambda $) and with the Haar-density topology of $\mathbb{R}_{+}.$ In
particular, the two topologies make available as a tool the \textit{%
interior-point theorem} of \textit{Steinhaus-Weil} from measure theory [Ste]
[Wei], and the \textit{Piccard-Pettis} category analogue [Pic] [Pet] (cf.
[BinO13]). Before identifying the (normalized) Haar measure of $\mathbb{G=G}%
_{\eta }$, written $\eta _{\mathbb{G}},$ we observe below that $\mathbb{G}$
has a natural order which coincides with the usual order on $\mathbb{R}.$ We
also identify the associated canonical invariant metrics on $\mathbb{G}$,
below.

We recall that by the Birkhoff-Kakutani Theorem ([Bir], [Kak1]; cf. [DieS,
\S 3.3, \S\ 8]) we may equip any metrizable group $G$ with a
(left-)invariant metric $d_{G}^{L}$, equivalently with a (group) \textit{norm%
} $||g||:=d_{G}^{L}(g,1_{G}),$ as in [BinO2] (`pre-norm' in [ArhT]) that
generates its topology. Its defining features are:

\noindent (i) $||g||=0$ iff $g=1_{G};$

\noindent (ii) $||gh||\leq ||g||+||h||;$

\noindent (iii) $||g^{-1}||=||g||.$

\bigskip

The group norm on $\mathbb{R}_{+}$ is also a limit of $||t||_{\rho }$ for $%
\rho \rightarrow \infty ,$ as we will see below.

\bigskip

\noindent \textbf{Proposition 1. }(a) \textit{A group-norm on }$\mathbb{G}%
_{\rho }$ \textit{for }$\rho \geq 0$ \textit{is given by}%
\[
||t||_{\rho }:=|\log (1+\rho t)|(1+\rho )/\rho .
\]%
(b) \textit{In particular, }$||1||_{\rho }=\log (1+\rho )/\rho ,$\textit{\
and }$||t||_{\rho }\rightarrow |t|$\textit{\ as }$\rho \rightarrow 0$\textit{%
\ (for }$t\neq 0).$

\noindent (c) \textit{A group-norm on }$\mathbb{G}_{\infty }=\mathbb{R}_{+}$
\textit{is given by}%
\[
||t||_{\infty }:=|\log t|.
\]

\bigskip

\noindent \textbf{Proof. }(a)\textbf{\ }Here (i) is clear; as for (ii), we
have
\[
||s\circ _{\rho }t||_{\rho }=|\log (1+\rho (s+t+\rho st))|(1+\rho )/\rho
=|\log (1+\rho s)(1+\rho t)|(1+\rho )/\rho \leq ||s||_{\rho }+||t||_{\rho }.
\]%
Then (iii) follows, since $\eta _{\rho }(t_{\rho }^{-1})=\eta _{\rho
}(t)^{-1}$ i.e. with $s=t_{\rho }^{-1}$%
\begin{equation}
(1+\rho s)=1/(1+\rho t).  \tag{inv}
\end{equation}%
(Or, from (conj) above, with $t$ for $y$ and its inverse $s=t_{\rho }^{-1}$
for $x,$
\[
1_{\rho }=0=s\circ _{\rho }t=[(1+\rho s)(1+\rho t)-1]/\rho :\qquad (1+\rho
s)(1+\rho t)=1.)
\]

(b)\textbf{\ }The second assertion follows by L'Hospital's rule (or as $\log
(1+\rho t)\sim \rho t$ for $\rho \sim 0)$.

(c) The final assertion is similar to but simpler than in (a). $\square $

\bigskip

\noindent \textbf{Remarks. }The inclusion of the scaling factor $(1+\rho )$
is dictated by Haar-measure normalization concerns, below.

\bigskip

\noindent \textbf{Proposition 2. }\textit{For }$\rho \geq 0,$ \textit{the
set }$\mathbb{G}_{\rho }^{+}=[0,\infty )$\textit{\ is a sub-semigroup of }$%
\mathbb{G}_{\rho };$\textit{\ the induced order, }$y\leq _{\rho }x$\textit{\
iff }$x\circ _{\rho }y^{-1}\in \lbrack 0,\infty ),$\textit{\ coincides with }%
$y\leq x.$ \textit{Furthermore, for }$c>0$\textit{\ and }$a<b,$%
\[
a\circ _{\rho }c\leq b\circ _{\rho }c;
\]%
\textit{in particular, for the interval }$(a,b)$,%
\[
(a,b)\circ _{\rho }c=(a\circ _{\rho }c,b\circ _{\rho }c):
\]%
\textit{the Euclidean topology on }$\mathbb{G}_{\rho }$\textit{\ is
invariant under (positive) translation by }$\circ _{\rho }.$

\noindent \textit{Likewise, for }$\rho >0,$ \textit{if }$a<b$\textit{\ and }$%
c<d,$\textit{\ with }$a,b,c,d\in \mathbb{G}_{\rho },$%
\[
a\circ _{\rho }c\leq b\circ _{\rho }d.
\]%
\textit{and}
\[
s\leq t\text{ iff }s_{\rho }^{-1}\geq t_{\rho }^{-1}\qquad (s,t\in \mathbb{G}%
_{\rho }).
\]

\bigskip

\noindent \textbf{Proof. }For the first assertion observe that
\[
0\leq x-(1+\rho x)y/(1+\rho y)\text{ iff }0\leq x(1+\rho y)-(1+\rho x)y=x-y,%
\text{ as }1+\rho y>0.
\]%
From here, as $a\leq b$ and $c\leq d,$%
\[
a\circ _{\rho }c\leq b\circ _{\rho }c\text{ and }c\circ _{\rho }b<d\circ
_{\rho }c:\text{\qquad }a\circ _{\rho }c\leq b\circ _{\rho }d.
\]%
Finally, $s\leq t$ iff
\[
-1/(1+\rho t)\geq -1/(1+\rho s):\text{\qquad }1-1/(1+\rho t)\geq 1-1/(1+\rho
s):\text{\quad }-\rho t_{\rho }^{-1}\geq -\rho s_{\rho }^{-1}.\text{\quad }%
\square
\]

\bigskip

\noindent \textbf{Theorem 1 (Haar measure). }\textit{Normalized Haar measure
on the Popa group }$\mathbb{G=G}_{\rho },$\textit{\ with the Euclidean
topology giving the interval }$(0,1)$\textit{\ measure }$||1||_{\rho }$
\textit{for }$\rho \geq 0,$ \textit{has Radon-Nikodym derivative }$(1+\rho
)/\eta _{\rho }(g)$\textit{\ w.r.t. }$\mathrm{d}g$\textit{, the Lebesgue
measure on the additive reals, that is}
\begin{eqnarray*}
\mathrm{d}\eta _{\mathbb{G}}(t) &=&(1+\rho )\text{ }\mathrm{d}t/\eta
(t)=(1+\rho )\text{ }\mathrm{d}t/\eta _{\rho }(t) \\
&=&(1+\rho )\text{ }\mathrm{d}t/(1+\rho t),\text{ for }\eta =\eta _{\rho }.
\end{eqnarray*}%
\textit{In particular, as }$1_{\rho }=0,$ \textit{the group norm satisfies}%
\[
||x||_{\rho }=\eta _{\mathbb{G}}((1_{\rho },x))=\int_{0}^{x}(1+\rho )\text{ }%
\mathrm{d}t/(1+\rho t)=\frac{1+\rho }{\rho }|\log (1+\rho x)|.
\]

\noindent \textbf{Proof. }Since Haar measure is unique up to
proportionality, begin by letting $\tilde{\eta}_{\mathbb{G}}$ be an
arbitrary Haar measure for the group. As $\tilde{\eta}_{\mathbb{G}}$ and $%
\lambda $ are absolutely continuous measures w.r.t. each other (both give
(non-degenerate) intervals positive measure), the Radon-Nikodym derivative,
which we write below as
\[
\delta (g):=\mathrm{d}\tilde{\eta}_{\mathbb{G}}/\mathrm{d}\lambda (g),
\]%
is well defined. To find the Radon-Nikodym derivative at $g,$ we compare the
Lebesgue measure of an interval around $g$ with its $\tilde{\eta}_{\mathbb{G}%
}$-measure. Taking $(a,b)$ an arbitrary interval around $0=1_{\rho },$ so
that $g\circ _{\eta }(a,b)$ is a neighbourhood of $g,$%
\[
g\circ _{\eta }(a,b)=g+\eta _{\rho }(g)(a,b):\qquad \lambda (g\circ _{\eta
}(a,b))=\eta _{\rho }(g)\lambda (a,b).
\]%
Now, taking limits below as $a\uparrow 0,b\downarrow 0,$ and setting $%
t=g\circ s=g+\eta (g)s$%
\begin{eqnarray*}
\delta (1_{\rho }) &=&\lim \frac{\tilde{\eta}_{\mathbb{G}}((a,b))}{\lambda
(a,b)}=\lim \frac{\tilde{\eta}_{\mathbb{G}}(g\circ _{\rho }(a,b))}{\lambda
(a,b)}\qquad \text{(invariance)} \\
&=&\lim \frac{\int_{g\circ _{\eta }(a,b)}\delta (t)\text{ }\mathrm{d}t}{%
\lambda (a,b)}=\eta _{\rho }(g)\lim \frac{\int_{(a,b)}\delta (g+\eta (g)s)%
\text{ }\mathrm{d}s}{\lambda (a,b)}=\eta _{\rho }(g)\delta (g)\qquad \text{%
a.e. }
\end{eqnarray*}%
by the Lebesgue differentiation theorem [Sak, IV \S\ 5], [Rud2, Th. 8.6]. So%
\[
\mathrm{d}\tilde{\eta}_{\mathbb{G}}(t)/\mathrm{d}t=\delta (g)=\delta
(1_{\rho })/\eta _{\rho }(g).\text{ }
\]%
So for the normalized measure $\eta _{\mathbb{G}}$ of the theorem, the
Radon-Nikodym derivative at $g$ is proportional to $1/\eta _{\rho }(g).$ The
proportionality constant $(1+\rho )$ allows for the two extreme $\rho $
values, to yield Lebesgue measure $\mathrm{d}t$ on the additive reals for $%
\rho =0,$ and Haar measure $\mathrm{d}t/t$ on the multiplicative reals $%
\mathbb{R}_{+}$ as $\rho \rightarrow \infty :$
\[
\mathrm{d}\eta _{\mathbb{G}}(t)=\frac{1+\rho }{1+\rho t}\text{ }\mathrm{d}%
t\rightarrow \mathrm{d}t\text{ as }\rho \rightarrow 0,\qquad \mathrm{d}\eta
_{\mathbb{G}}(t)=\frac{1+\rho }{1+\rho t}\text{ }\mathrm{d}t\rightarrow
\frac{\mathrm{d}t}{t}\text{ as }\rho \rightarrow \infty .\text{ }\square
\]

\bigskip

\noindent \textbf{Remark. }For $\eta =\eta _{\rho }$ and $\rho =0,$ we
interpret $\rho ^{\ast }=-1/\rho $ to mean $-\infty $ (the Popa centre
recedes to $-\infty $); then, since $s\circ _{0}t=s+t,$ we recover the
additive reals under ordinary Lebesgue measure, so $\mathbb{G}_{0}=\mathbb{R}
$, by Prop. 1. Here, computing distance relative to $1_{\rho }=0,$
\[
||x||_{\rho }=\eta _{\rho }(0,x)=\int_{0}^{x}\frac{1+\rho }{1+\rho t}dt=%
\frac{1+\rho }{\rho }|\log (1+\rho x)|\rightarrow |x|\text{ as }\rho
\rightarrow 0
\]%
(the modulus signs are needed iff $x<0,$ when $(x,0)$ replaces $(0,x)$
above). As before, $\log (1+\rho x)\sim \rho x$ for $\rho \sim 0.$

In the limit as $\rho \rightarrow \infty $ we interpret $\rho ^{\ast
}=-1/\rho $ to mean $0$ (the Popa centre approaches $0,$ from the left).
Since $1_{\rho =\infty }=1$ is the unit in the multiplicative reals $\mathbb{%
R}_{+}:=(0,\infty )$, computing distance now relative to $1$, we retrieve
\begin{eqnarray*}
\eta _{\rho }(1,x) &=&\int_{1}^{x}\frac{1+\rho }{1+\rho t}\text{ }\mathrm{d}%
t=\frac{1+\rho }{\rho }\log \left\vert \frac{1+\rho x}{1+\rho }\right\vert
\rightarrow |\log x|\text{ as }\rho \rightarrow \infty : \\
||x||_{\rho =\infty } &=&|\log x|.
\end{eqnarray*}%
Recalling that
\[
x\circ _{\rho }y=x+y(1+\rho x)=x+y+\rho xy=[(1+\rho x)(1+\rho y)-1]/\rho ,
\]%
the corresponding conjugacy yields%
\[
(x\circ _{\rho }y)/\rho =[(1+\rho x)(1+\rho x)-1]/\rho ^{2}=x/\rho +y/\rho
+xy\rightarrow xy,\text{ as }\rho \rightarrow \infty ,
\]%
so that $\mathbb{G}_{\infty }$ has domain $\mathbb{R}_{+}$ with $\circ
=\times $ (`the multiplicative reals').

This means that, up to scaling, there are just \textit{three} Popa
operations/groups, corresponding to $\rho =0,1,\infty ,$ namely $+,\circ
,\times $ with $\circ $ the circle operation of ring theory as above.

\bigskip

\noindent \textbf{Remarks. }1. The alternative normalization is $\delta
(1_{\rho })=1,$ as%
\[
x\circ _{\rho }y=x+y+\rho xy\sim x+y\quad (x,y\rightarrow 0).
\]

\noindent \textbf{2.} Note that $\mathbb{G}_{\rho }$ for $\rho \in \mathbb{R}%
_{+}$ has only one idempotent, $c=0$ (replaced by $c=1$ in the case $\rho
=\infty ):$
\begin{eqnarray*}
c &=&c\circ _{\rho }c:\qquad 0=c+\rho c^{2}=c(1+\rho c), \\
c &=&c^{2}:\qquad c=1\text{ in }\mathbb{R}_{+}=(0,\infty ).
\end{eqnarray*}

\noindent \textbf{3.} The origin of the Haar measures $\mathrm{d}t$, $%
\mathrm{d}t/t$ for the cases $\rho =0,\infty $ above are clear: the
arithmetic operations $+$ and $\times $. From the canonical $\mathrm{d}t/t$
case one may infer the general $\rho \in (0,\infty )$ case by a change of
origin to the Popa centre $-1/\rho .$ That of the intermediate values $\rho
\in (0,\infty )$ is exemplified by the case $\rho =1$, giving $\mathrm{d}%
t/(1+t)$. This arises via the role of the Beck sequences in the proof of
Theorem 3 and the Remark below it, and is an instance of the ergodic theorem
(see e.g. Billingsley [Bil, Ch. 1 \S 4] and Remark 4 below). The same
measure arises in the Gauss-Kuzmin theorem on continued fractions, and for
the same reason (again, see [Bil, Ch. 1 \S 4]).

\noindent \textbf{4.} As mentioned above: see [BinO7, Prop. 11 (iv)] for the
sense in which the Beck sequence of iterates above grows arithmetically,
which links their averages with the arithmetic means in the
(Birkhoff-Khinchin) ergodic theorem.

\noindent \textbf{5.} The limiting behaviour of the moving average $%
[U(x\circ _{\varphi }t)-U(x)]/\varphi (x)$ of $U$ and the Tauberian
one-sided conditions studied by Bingham and Goldie [BinG] emerge in \S 3
below directly from the asymptotic operator $K_{h}^{\varphi }(t,x)$ with the
specialization $h=\varphi .$ The group norms exhibited in Th. 1 above thus
coincide with the measures of occupation `times' (on $[1_{\rho },x]$) of the
associated limiting velocity flow $\mathrm{d}w(t)/\mathrm{d}t=\eta (t)$ for $%
\eta =\eta ^{\varphi }.$ Here Lebesgue measure $\mathrm{d}t$ measures time,
and equates with $w^{\prime }(t)\cdot \mathrm{d}t/\eta (t),$ i.e. the Haar
integral of the flow rate.

\bigskip

We recall that the dual of a locally compact abelian group $G$, denoted $%
\hat{G},$ comprises the continuous homomorphisms from $G$ to $\mathbb{T}$,
the unit circle in the complex plane $\mathbb{C}$. For $\eta _{G}$ a Haar
measure on $G,$ the Fourier transform is defined by
\[
\hat{f}(\gamma ):=\int_{G}f(g)\gamma (-g)\text{ }\mathrm{d}\eta
_{G}(g)\qquad (\gamma \in \hat{G});
\]%
for background see [Rud1], [Loo]. We specialize this in Theorem 2 below to
the Popa group $(\mathbb{G}_{\rho },\circ _{\rho })$ for $0<\rho <\infty .$
It is helpful to first consider the extreme cases $\rho =0$ and $\rho
=\infty ,$ corresponding respectively to the familiar cases $G=(\mathbb{R}%
,+) $ and $G=(\mathbb{R}_{+},\times ).$ In the first case $\hat{G}=G=(%
\mathbb{R},+)$ [Loo, 35C], and we may write
\[
\gamma (w)=e^{i\gamma w}\qquad (\gamma \in \mathbb{R}),
\]%
so that, for $f\in L^{1}(\mathbb{R}),$%
\[
\hat{f}(\gamma )=\int_{\mathbb{R}}f(w)e^{-i\gamma w}\text{ }\mathrm{d}%
w\qquad (\gamma \in \mathbb{R}).
\]%
We pass to the second case using the isomorphism $w=\log v$ which, for $f\in
L^{1}(\mathbb{R}_{+}),$ yields both the Fourier and Mellin transforms as%
\[
\hat{f}(\gamma )=\int_{0}^{\infty }f(v)e^{-i\gamma \log v}\text{ }\mathrm{d}%
v/v\quad (\gamma \in \mathbb{R}),\qquad \check{f}(z)=\int_{0}^{\infty
}f(t)t^{-z}\text{ }\mathrm{d}t/t\quad (z\in \mathbb{C}),
\]%
with characters represented multiplicatively by $\gamma (t)=t^{z}.$

We turn to the Fourier transform in the context of a locally compact abelian
group ([Rud1], [Loo]), specialized to the Popa-group $\mathbb{G}_{\rho }$
for $\rho >0.$ As we shall see, the Fourier-Popa transform of $f:\mathbb{G}%
_{\rho }\rightarrow \mathbb{R}$ is in fact the ordinary Fourier transform of
an affinely related function $f_{\rho }:\mathbb{R}_{+}\rightarrow \mathbb{R}$%
, defined as follows:
\[
f_{\rho }(t)=\frac{1+\rho }{\rho }f(\eta _{\rho }^{-1}(t))=\frac{1+\rho }{%
\rho }f((t-1)/\rho ),
\]%
so that%
\[
f_{\rho }(1/u)=\frac{1+\rho }{\rho }f(\eta _{\rho }^{-1}(u)_{\rho }^{-1})=%
\frac{1+\rho }{\rho }f((1-u)/(\rho u)).
\]%
As expected, for $\rho \rightarrow \infty $ we recover $f$ by rescaling: $%
f_{\rho }(\rho t)\rightarrow f(t)$.

\bigskip

\noindent \textbf{Theorem 2 (Fourier transform).} \textit{For the Popa group}
$G=(\mathbb{G}_{\rho },\circ _{\rho })$ \textit{with }$0<\rho <\infty ,$%
\textit{\ the characters }$\gamma \in \hat{G}$\textit{\ are }%
\[
\gamma (u):=e^{i\gamma \log (1+\rho u)}\qquad (\gamma \in \mathbb{R}).
\]%
\textit{So, writing }$+_{\rho }$\textit{\ and }$-_{\rho }$ \textit{for the
operations of }$\circ _{\rho }$\textit{\ and inversion here, the Fourier
transform corresponding to the canonical Haar measure of Theorem 1 is }%
\[
\hat{f}(\gamma )=\int_{\mathbb{G}_{\rho }}f(u)\gamma (-_{\rho }u)(1+\rho )%
\text{ }\mathrm{d}u/(1+\rho u)=\int_{0}^{\infty }f_{\rho }(t)e^{i\log
t^{-\gamma }}\text{ }\mathrm{d}t/t,
\]%
\textit{that is}
\[
\hat{f}(\gamma )=\int_{\mathbb{G}_{\rho }}f(u)\gamma (-_{\rho }u)(1+\rho )%
\text{ }\mathrm{d}u/(1+\rho u)=\int_{0}^{\infty }f((t-1)/\rho )e^{i\log
t^{-\gamma }}\text{ }\mathrm{d}t/t.
\]%
\textit{The corresponding Mellin transform is thus}%
\[
\check{f}(z)=\int_{0}^{\infty }f_{\rho }(t)t^{-z}\text{ }\mathrm{d}%
t/t=\int_{0}^{\infty }f_{\rho }(1/u)u^{z}\text{ }\mathrm{d}%
u/u=\int_{0}^{\infty }f((1-u)/(\rho u))u^{z}\text{ }\mathrm{d}u/u.
\]

\bigskip

\noindent \textbf{Proof. }Applying the isomorphisms $\eta _{\rho }:(\mathbb{G%
}_{\rho },\circ )\rightarrow (\mathbb{R}_{+},\times )$ and $\log :(\mathbb{R}%
_{+},\times )\rightarrow (\mathbb{R},+)$ and using $u,v,w$ as corresponding
generic elements with $w=\log v$ and $v=1+\rho u,$ the character
representation for $(\mathbb{R},+)$ recalled above gives the character
representation for $(\mathbb{G}_{\rho },\circ )$ as asserted. By (inv) above
\[
1+\rho (-_{\rho }u)=1/(1+\rho u),
\]%
so substitution for $\gamma (-_{\rho }u)$ gives the Fourier transform as
\[
\hat{f}(\gamma )=\int_{\mathbb{G}_{\rho }}f(u)\gamma (-_{\rho }u)(1+\rho )%
\text{ }\mathrm{d}u/(1+\rho u)=\int_{-1/\rho }^{\infty }f(u)e^{-i\gamma \log
(1+\rho u)}(1+\rho )\text{ }\mathrm{d}u/(1+\rho u).
\]%
Putting $t=\eta _{\rho }(u)=1+\rho u$ gives
\[
\hat{f}(\gamma )=\int_{0}^{\infty }f_{\rho }(t)e^{i\log t^{-\gamma }}\text{ }%
\mathrm{d}t/t.
\]%
This gives the first form of the Mellin transform above; for the second,
take $u=1/t$. $\square $

\section{Asymptotic actions and functional equations}

We begin with the \textit{Karamata asymptotic operator }$K$ acting on $f:%
\mathbb{R}_{+}\rightarrow \mathbb{R}_{+},$ as in ($K$) of \S\ 1:%
\[
K(t,x)f:=\frac{f(xt)}{f(x)}.
\]%
Suppose that $f$ is \textit{Karamata regularly varying}, i.e. that, as%
\textit{\ }$x\rightarrow \infty ,$\textit{\ }%
\[
K(t,x)f:=\frac{f(xt)}{f(x)}\rightarrow K_{f}(t).
\]%
Here we adopt a relatively new point of view on the classical theory, by
making explicit use of what has so far been mostly implicit: the\textit{\
cocycle structure} underlying the operator\textit{\ }$K(t,x),$ cf. [Ell]
[EllE]. It is this that characterizes the limit function $K_{f},$ the
\textit{Karamata kernel }of\textit{\ }$f$. Indeed,%
\[
\frac{f(xts)}{f(x)}=\frac{f(xts)}{f(xt)}\cdot \frac{f(xt)}{f(x)}:\qquad
K(st,x)=K(s,xt)K(t,x).
\]%
In the limit this yields the multiplicative \textit{Cauchy functional
equation},%
\begin{equation}
K_{f}(st)=K_{f}(s)K_{f}(t).  \tag{$CFE$}
\end{equation}%
We will need the \textit{Popa operation} $\circ _{h}$ above to be
associative, and (see Th. O below) this requires $h$ to satisfy the\textit{\
Go\l \k{a}b-Schintzel equation}:%
\begin{equation}
h(s\circ _{h}t)=h(s+h(s)t)=h(s)h(t).  \tag{$GS$}
\end{equation}%
Thus ($GS$) expresses homomorphy in this context, which will occur after the
passage to the limit $x\rightarrow \infty $. Before taking this limit, one
has instead `asymptotic associativity', or `almost associativity'. The Popa
notation $x\circ _{\varphi }t=x+t\varphi (x)$ describes a $t$-translation
modified locally at $x,$ or `accelerated at $x$' by reference to the
`accelerator' $\varphi $\textit{\ (positive)}. We will need the rate of
acceleration and its asymptotic value for the $t$-translation:%
\[
\eta _{x}(t),\text{ or }\eta _{x}^{\varphi }(t),:=\frac{\varphi (x\circ
_{\varphi }t)}{\varphi (x)}=\frac{\varphi (x+t\varphi (x))}{\varphi (x)}%
\rightarrow \eta (t),\text{ or }\eta ^{\varphi }(t)
\]%
(assumed to exist), so that $\eta (t)\geq 0.$ As we learn from the Uniform
Convergence Theorem (UCT) below, for $\varphi $ above Baire or measurable,
convergence is necessarily \textit{locally uniform}. The relevance of such
convergence is witnessed by

\bigskip

\noindent \textbf{Theorem O }[Ost1, Th. 0].\textbf{\ }\textit{If }$\varphi
(x)=O(x)$\textit{\ and }$\eta _{x}(t)\rightarrow \eta (t)=\eta ^{\varphi
}(t),$ \textit{locally uniformly in }$t$\textit{, then }$\eta $\textit{\
satisfies the Go\l \k{a}b-Schinzel functional equation}%
\begin{equation}
\eta (s\circ _{\eta }t)=\eta (s)\eta (t).  \tag{$GS$}
\end{equation}

\bigskip

\noindent \textit{Notational conventions. }In Theorem O above $\eta _{x}$
contains the $x$ which tends to infinity. After this passage to the limit,
attention focuses on the limit function $\eta (t)$ which will depend on a
parameter $\rho $, below. We allow ourselves to denote this limit by $\eta
_{\rho }(t)$ and let context speak for itself here. Below we will take $%
GS:=\{\eta _{\rho }:\rho \geq 0\}$ to denote the family of continuous
(positive) solutions of the equation ($GS$).

\bigskip

For $\varphi $ Baire/measurable $\eta ^{\varphi }$ is likewise
Baire/measurable and so, as a solution of ($GS$), \textit{continuous}, by a
theorem of Popa [Pop]. Furthermore, non-negative solutions of ($GS$), being
bounded below, are likewise continuous, as noted in \S 2. In any case, here
we are interested only in \textit{positive} solutions of ($GS$), and these
take the form $\eta (t)=\eta _{\rho }(t):=1+\rho t,$ for $t>\rho ^{\ast
}:=-1/\rho $ with $\rho \geq 0$ (and $0$ to the left of $\rho ^{\ast },$
though here we work in $\mathbb{R}_{+}$), by a theorem of Go\l \k{a}b and
Schinzel -- for the literature see [Brz1], [Jab5], and [Ost1]. For a
discussion of circumstances when local boundedness implies the continuity of
solutions, for the family relevant here of functional equations related to ($%
GS$), see [Jab3].

Below, we will encounter \textit{two }auxiliary functions, $h$ and $\varphi
, $ the second of which will give such an $\eta $ asymptotically (so $\eta $
satisfies ($GS$) and $\circ _{\eta }$ is associative).

For the purposes of combining an $s$- and a $t$-translation, it is
convenient to expand the accelerator notation to one parametrized locally at
$x$:%
\[
s\circ _{\varphi x}t:=s+t\eta _{x}(s)=s+t\frac{\varphi (x+s\varphi (x))}{%
\varphi (x)}.
\]%
So in the limit one has for $\eta =\eta ^{\varphi }=\eta _{\rho }$:%
\[
\circ _{\varphi x}\rightarrow \circ _{\eta }=\circ _{\rho }.
\]%
This justifies our earlier reference to `asymptotic associativity'. A second
reason for the term comes from a very convenient expression for a related
form of associativity, one which otherwise the notation keeps hidden:%
\[
(x\circ _{\varphi }b)\circ _{\varphi }a=x\circ _{\varphi }(b\circ _{\varphi
x}a)
\]%
As an immediate application of this framework, we can rephrase the Beurling
asymptotics, clarifying the underlying \textit{cocycle structure}. These, as
we will see, lead to functional equations, whose solutions are discussed in
\S 5 below -- see also the surveys [Brz1] and [Jab5]; cf. [Ost1].

\bigskip

\noindent \textbf{Proposition 3 (Beurling regular variation). }\textit{For
the Beurling asymptotic operator }$K^{\varphi }$ \textit{acting on} $f:%
\mathbb{R}_{+}\rightarrow \mathbb{R}_{+},$\textit{\ }%
\[
K^{\varphi }(t,x)f:=\frac{f(x+t\varphi (x))}{f(x)}=\frac{f(x\circ _{\varphi
}t)}{f(x)},
\]%
\textit{suppose that }$f$\textit{\ is Beurling regularly varying, i.e. that,
as }$x\rightarrow \infty ,$\textit{\ }%
\[
K^{\varphi }(t,x)f:=\frac{f(x+t\varphi (x))}{f(x)}=\frac{f(x\circ _{\varphi
}t)}{f(x)}\rightarrow K_{f}(t).
\]%
\textit{The corresponding cocycle structure is}%
\[
K^{\varphi }(t\circ _{\varphi x}s,x)=K^{\varphi }(s,x\circ _{\varphi
}t)K^{\varphi }(t,x),
\]%
\textit{leading in the limit to the Chudziak-Jab\l o\'{n}ska equation}%
\begin{equation}
K_{f}(t\circ _{\eta }s)=K_{f}(s)K_{f}(t).  \tag{$CJ$}
\end{equation}

\noindent \textbf{Proof. }We have
\[
\frac{f(x+(s+t)\varphi (x))}{f(x)}=\frac{f(x+t\varphi (x)+(s/\eta _{x})\cdot
\varphi (x+t\varphi (x)))}{f(x+t\varphi (x))}\cdot \frac{f(x+t\varphi (x))}{%
f(x)},
\]%
so that in the limit%
\[
K(s+t,x)=K(s/\eta _{x}(t),x+t\varphi (x))K(t,x).
\]%
Here replacing $s$ by $s\eta _{x}(t)$ yields%
\[
K(t+s\eta ^{x}(t),x)=K(s,x+t\varphi (x))K(t,x).\qquad \square
\]

We turn now to the \textit{general regular variation }of the title and \S 1
(cf. [BinO14]).

Following [Ost1], the auxiliary function $\varphi :\mathbb{R}_{+}\rightarrow
\mathbb{R}_{+}$ is \textit{self-equivarying}, $\varphi \in SE,$ if $\varphi
(x)=O(x)$ and $\eta _{x}^{\varphi }(t)\rightarrow \eta (t)=\eta ^{\varphi
}(t),$ locally uniformly in\textit{\ }$t,$ as in Theorem O. The auxiliary
function $h$ will be Beurling regularly varying as in Prop. 1, i.e. $\varphi
$-regularly varying, in the sense of [BinO5].

\bigskip

\noindent \textbf{Proposition 4 (General regular variation). }\textit{For
the general asymptotics\ }%
\[
K_{h}^{\varphi }(t,x):=\frac{f(x+t\varphi (x))-f(x)}{h(x)}\rightarrow
K_{f}(t),
\]%
\textit{with }$\varphi \in SE,$ \textit{the corresponding cocycle structure
is}
\[
K_{h}^{\varphi }(t+s\eta _{x}(t),x)=K_{h}^{\varphi }(t\circ _{\varphi
x}s,x)=K_{h}^{\varphi }(s,x)K^{\varphi }(t,x)+K_{h}^{\varphi }(t,x),
\]%
\textit{leading in the limit to}%
\[
K_{f}(t+s\eta (t))=K_{f}(s)K_{h}(t)+K_{f}(t),
\]%
\textit{or, equivalently, to the Beurling-Goldie equation satisfied by }$%
K_{f}:\mathbb{G}_{\eta }\rightarrow \mathbb{G}_{\sigma }:$%
\begin{equation}
K_{f}(t\circ _{\eta }s)=K_{f}(t)\circ _{\sigma }K_{f}(s),\qquad \text{for }%
\sigma (z)=K_{h}(K_{f}^{-1}(z)).  \tag{$BG$}
\end{equation}

\bigskip

\noindent \textbf{Proof. }Here the underlying cocycle structure mixes
products with addition: with $y:=x\circ _{\varphi }t,$%
\begin{eqnarray*}
K_{h}^{\varphi }(s+t,x) &=&\frac{f(x+(s+t)\varphi (x))-f(x)}{h(x)} \\
&=&\frac{f(x+t\varphi (x)+(s/\eta _{x})\varphi (x\circ _{\varphi
}t))-f(x\circ _{\varphi }t)}{h(x\circ _{\varphi }t)}\frac{h(x\circ _{\varphi
}t)}{h(x)}+\frac{f(x\circ _{\varphi }t)-f(x)}{h(x)} \\
&=&\frac{f(y+(s/\eta _{x})\varphi (y))-f(y)}{h(y)}\frac{h(x\circ _{\varphi
}t)}{h(x)}+K_{h}^{\varphi }(t,x) \\
&=&K_{h}^{\varphi }(s/\eta _{x},y)K^{\varphi }(t,x)+K_{h}^{\varphi }(t,x).
\end{eqnarray*}%
In the limit, since $x+t\varphi (x)=x(1+t\varphi (x)/x)\rightarrow \infty $
and $\varphi (x)=O(x),$

\[
K_{f}(s+t)=K_{f}(s/\eta ,x)K_{h}(t)+K_{f}(t),
\]%
giving $(BG)$ as above. $\square $

\bigskip

\noindent \textbf{Remark. }A measurable $\varphi :\mathbb{R}_{+}\rightarrow
\mathbb{R}_{+}$ is said to be \textit{Beurling slowly varying} if, as above,
but with $\varphi (x)=o(x)$ and $\eta ^{\varphi }(t)\equiv 1$ (that is, $%
\rho =0$ in the above); it is \textit{self-neglecting} if the convergence $%
\eta _{x}(t)\rightarrow 1$ is locally uniformly in\textit{\ }$t$ -- see
[BinGT \S\ 2.11], [BinO5].

\section{Subadditivity in Popa groups}

\noindent \textbf{Definition. }For\textbf{\ }$\rho ,\sigma \in \lbrack
0,\infty ],$ call $S:\mathbb{G}_{\rho }\rightarrow \mathbb{G}_{\sigma }$
\textit{subadditive} (resp. \textit{additive}) if
\[
S(x\circ _{\rho }y)\leq S(x)\circ _{\sigma }S(y)\qquad \text{resp.}\qquad
S(x\circ _{\rho }y)=S(x)\circ _{\sigma }S(y),
\]%
or in the notation of Theorem 2%
\[
S(x+_{\rho }y)\leq S(x)+_{\sigma }S(y)\qquad \text{resp.}\qquad S(x+_{\rho
}y)=S(x)+_{\sigma }S(y),
\]

As $\mathbb{G}_{0}=\mathbb{R}$ (the additive reals), when $\rho =\sigma =0,$
this yields the usual notion of subadditivity, resp. additivity.

In particular the solutions $K:\mathbb{G}_{\rho }\rightarrow \mathbb{G}%
_{\sigma }$ to the equation $(BG)$ are additive. For fixed $\rho ,\sigma \in
\mathbb{R}_{+}$ with $\sigma >0,$ the canonical form depends on a parameter $%
\kappa \in \mathbb{R}$ (Theorem 3 below, [Ost2-Hom], [Chu1,2]), as follows:%
\begin{eqnarray*}
K_{\kappa }(t) &=&\eta _{\sigma }^{-1}(\eta _{\rho }(t)^{\kappa }) \\
&=&[(1+\rho t)^{\kappa }-1]/\sigma ,\text{ if also }\rho >0.
\end{eqnarray*}%
Above one has $\eta _{\rho }:\mathbb{G}_{\rho }\rightarrow \mathbb{R}_{+}$,
and $\eta _{\sigma }^{-1}:\mathbb{R}_{+}\rightarrow \mathbb{G}_{\sigma }.$
The case $\kappa =0$ corresponds to the trivial solution $K\equiv 1_{\sigma
}=0.$

\bigskip

\noindent \textbf{Example}. Recalling that $\eta _{\rho }(x\circ _{\rho
}y)=\eta _{\rho }(x)\eta _{\rho }(y),$ so that $\eta _{\rho }^{-1}(uv)=\eta
_{\rho }^{-1}(u)\circ _{\rho }\eta _{\rho }^{-1}(v)$ (on substituting $%
u=\eta _{\rho }(x)$ etc.),
\begin{eqnarray*}
K(x\circ _{\rho }y) &=&\eta _{\sigma }^{-1}(\eta _{\rho }(x\circ _{\rho
}y)^{\kappa })=\eta _{\sigma }^{-1}(\eta _{\rho }(x)^{\kappa }\eta _{\rho
}(y)^{\kappa }) \\
&=&\eta _{\sigma }^{-1}(\eta _{\rho }(x)^{\kappa })\circ _{\sigma }\eta
_{\sigma }^{-1}(\eta _{\rho }(y)^{\kappa }) \\
&=&K(x)\circ _{\sigma }K(y).
\end{eqnarray*}

In fact, for fixed $\rho ,\sigma \in \mathbb{R}_{+},$ the only additive
functions bounded above are of this form, as below. Theorem 3 below is our
reformulation here of [Ost2-Hom, Prop. A]; cf. [Chu1,2].

\bigskip

\noindent \textbf{Theorem 3.} \textit{Take }$\psi :\mathbb{G}_{\rho
}\rightarrow \mathbb{G}_{\sigma }$\textit{\ additive with }$\rho ,\sigma \in
\lbrack 0,\infty ].$\textit{\ Then the lifting }$\Psi :\mathbb{R}\rightarrow
\mathbb{R}$\textit{\ defined by the canonical isomorphisms }$\log ,\exp ,$%
\textit{\ }$\{\eta _{\rho }:\rho >0\}$\textit{\ of }$\psi $\textit{\ to }$%
\mathbb{R}$\textit{\ is bounded above on }$\mathbb{G}_{\rho }$\textit{\ iff }%
$\Psi $\textit{\ is bounded above on }$\mathbb{R}$\textit{, in which case }$%
\Psi $ \textit{and }$\psi $\textit{\ are continuous. Then for some }$\kappa
\in \mathbb{R}$ \textit{one has:}

\renewcommand{\arraystretch}{1.25}%
\[
\begin{tabular}{|l|l|l|l|}
\hline
Popa parameter & $\sigma =0$ & $\sigma \in (0,\infty )$ & $\sigma =\infty $
\\ \hline
$\rho =0$ & $\kappa t$ & $\eta _{\sigma }^{-1}(e^{\kappa t})$ & $e^{\kappa
t} $ \\ \hline
$\rho \in (0,\infty )$ & $\log \eta _{\rho }(t)^{\kappa }$ & $\eta _{\sigma
}^{-1}(\eta _{\rho }(t)^{\kappa })$ & $\eta _{\rho }(t)^{\kappa }$ \\ \hline
$\rho =\infty $ & $\log t^{\kappa }$ & $\eta _{\sigma }^{-1}(t^{\kappa })$ &
$t^{\kappa }$ \\ \hline
\end{tabular}%
\]%
\medskip \newline
\renewcommand{\arraystretch}{1}

\noindent \textbf{Proof. }The canonical isomorphisms are order-preserving
and continuous. For $\rho ,\sigma >0$ the lifting is given by
\[
\Psi .=\log \eta _{\sigma }\psi \eta _{\rho }^{-1}\exp .,
\]%
and this still holds for extreme values of $\rho ,\sigma $ with $\exp ,$ $%
\log $ replacing $\eta _{0},\eta _{\infty }$. For $\Psi (x)=\kappa x,$ a
routine calculation gives $\psi $ as in the table above. $\square $

\bigskip

\noindent \textbf{Remark. }Notice that the passage from first to the third
column is effected via $\exp /\log ,$ while the middle column to the first
column requires scaling of the domain via the coefficient $\kappa :$

\begin{eqnarray*}
\lim\nolimits_{\sigma \rightarrow 0}\eta _{\sigma }^{-1}(e^{t\sigma \kappa
}) &=&\lim\nolimits_{\sigma \rightarrow 0}\frac{e^{\kappa \sigma t}-1}{%
\sigma }=\kappa t\text{ (equiv. }\kappa t\sim \log (1+\kappa \sigma
t)/\sigma \text{);} \\
\lim\nolimits_{\sigma \rightarrow 0}\eta _{\sigma }^{-1}(\eta _{\rho
}(t)^{\kappa \sigma }) &=&\lim\nolimits_{\sigma \rightarrow 0}\frac{[(1+\rho
t)^{\sigma \kappa }-1]}{\sigma }=\log \eta _{\rho }(t)^{\kappa }; \\
\lim\nolimits_{\sigma \rightarrow 0}\eta _{\sigma }^{-1}(t^{\kappa \sigma })
&=&\lim\nolimits_{\sigma \rightarrow 0}\frac{e^{\kappa \sigma \log t}-1}{%
\sigma }=\log t^{\kappa }.
\end{eqnarray*}

Note also%
\[
\text{ }\kappa t\sim \rho \log (1+\kappa t/\rho ),\text{ as }\rho
\rightarrow \infty \text{;}\qquad \lim\nolimits_{\rho \rightarrow 0}\frac{%
\log \eta _{\rho }(t)^{\kappa }}{\rho }=\lim\nolimits_{\rho \rightarrow 0}%
\frac{\kappa \log (1+\rho t)}{\rho }=\kappa t.
\]

\bigskip

\noindent \textbf{Definition. }Call $S:\mathbb{G}_{\rho }\rightarrow \mathbb{%
G}_{\sigma }$ \textit{additively bounded} on $\Sigma $ if for some $\kappa $%
\[
S(t)\leq K_{\kappa }(t)\qquad (t\in \Sigma ).
\]%
This lifts to the Popa context the notion of linear boundedness used in
[BinO10].

In the results below recall that $0=1_{\rho }=1_{\sigma };$ $B_{\delta }(x)$
is the open ball about $x$ of radius $\delta .$

\bigskip

\noindent \textbf{Proposition 5.} \textit{For} $S:\mathbb{G}_{\rho
}\rightarrow \mathbb{G}_{\sigma }$\textit{\ subadditive:\newline
\noindent }(i) \textit{if }$S$\textit{\ is bounded above on some interval,
say by }$M$\textit{\ on }$B_{\delta }(a)$\textit{, then for any }$b\in
\mathbb{G}_{\rho }^{+}$%
\[
S(b\circ a)\circ _{\sigma }M_{\sigma }^{-1}\leq S(x)\leq S(b\circ a_{\rho
}^{-1})\circ _{\sigma }M\qquad (x\in B_{\delta }(b))
\]%
\textit{\ (with }$M_{\sigma }^{-1}$\textit{\ etc. the inverses in the
corresponding groups); in particular it is locally bounded.\newline
}\noindent (ii) \textit{If }$S$\textit{\ is locally bounded, then }$\lim
\inf_{t\rightarrow 0}S(t)\geq 0,$ \textit{so }$S(0+)=0$\textit{\ if }$%
(HS(S)) $\textit{\ holds.}

\bigskip

Below (as in \S 1), `G for Goldie, G for general':

\bigskip

\noindent \textbf{Theorem G1. }\textit{For subadditive} $S:\mathbb{G}_{\rho
}^{+}\rightarrow \mathbb{G}_{\sigma }^{+}\cup \{-\infty ,+\infty \}$ \textit{%
with }$S(0+)=S(0)=0:$\textit{\ }$S$\textit{\ is continuous at }$0$ \textit{%
iff }$S(z_{n})\rightarrow 0,$\textit{\ for some sequence }$z_{n}\uparrow 0,$%
\textit{\ and then }$S$ \textit{is continuous everywhere, if finite-valued.}

\bigskip

\noindent \textbf{Proof of Theorem G1. }This is as in [BinO10], mutatis
mutandis, as the group order is the usual order on the line (Prop. 4), and
with $-x$ etc. replaced by $x_{\rho }^{-1}$ (equivalently by $-_{\rho }x$ as
in Theorem 2). It is critical here that one works in $\mathbb{G}_{\rho }^{+}$
and $\mathbb{G}_{\sigma }^{+}.$ $\square $

\bigskip

\noindent \textbf{Theorem G2\ }[BinO8, Th. 3]\textbf{.} \textit{If }$S:%
\mathbb{G}_{\rho }\rightarrow \mathbb{G}_{\sigma }$ \textit{\ is subadditive
with }$S(0)=0$ \textit{and there is a symmetric set }$\Sigma $\textit{\
containing }$0$ \textit{with:}

\noindent (i)\textit{\ }$S$\ \textit{is continuous at }$0$\textit{\ on }$%
\Sigma $\textit{;}

\noindent (ii) \textit{for all small enough }$\delta >0,$ $\Sigma
_{0}^{\delta }$ \textit{is locally Steinhaus-Weil}\newline
\noindent -- \textit{then }$S$ \textit{is continuous at }$0$\textit{\ and so
everywhere. }

\textit{In particular, this conclusion holds if there is a symmetric set }$%
\Sigma $\textit{, Baire/measurable and non-negligible in each }$(0,\delta )$%
\textit{\ for }$\delta >0,$\textit{\ on which}%
\[
S(u)=K_{\kappa _{\pm }}(u)\text{ }\mathit{for\ some\ }\kappa _{\pm }\in
\mathbb{R}\mathit{\ and\ all\ }u\in \mathbb{G}_{\rho }^{+}\cap \Sigma
\mathit{,or}\text{ }\mathit{all}\text{ }u\in \mathbb{G}_{\rho }^{-}{}\cap
\Sigma \text{ }\mathit{resp.}
\]

\bigskip

\noindent \textbf{Proof of Theorem G2.} W.l.o.g. $\sigma >0,$ as the case $%
\sigma =0$ is similar but simpler. Since $S|\Sigma $ is continuous at $0$ it
is bounded above on $\Sigma _{\delta }:=\Sigma \cap (\delta _{\sigma
}^{-1},\delta )$ for some $\delta >0;$ but $\Sigma _{\delta }\circ \Sigma
_{\delta }$ contains an interval, so $S$ is bounded on an interval, and so
locally bounded by Prop. 5(i). If $S$ is not continuous at $0$, then by
Prop. 5(ii) $\lambda _{+}:=\lim \sup_{t\rightarrow 0}S(t)>\lim
\inf_{t\rightarrow 0}S(t)\geq 0.$ Choose a null sequence $\{z_{n}\}$ with $%
S(z_{n})\rightarrow \lambda _{+}>0.$ Let $\varepsilon :=\min \{\lambda
_{+}/6,1/\sigma \}.$ W.l.o.g. $S(z_{n})>\lambda _{+}-\varepsilon $ for all $%
n.$ By continuity on $\Sigma $ at $0,$ there is $\delta >0$ with $%
|S(t)|<\varepsilon $ for $t\in \Sigma _{\delta }.$ As before and using
symmetry, $\Sigma _{\delta }\circ \Sigma _{\delta }=\Sigma _{\delta }\circ
(\Sigma _{\delta })_{\sigma }^{-1}$ contains an interval $I$ around $0.$ For
any $n$ with $z_{n}\in I,$ there are $u_{n},v_{n}\in \Sigma _{\delta }$ with
$z_{n}=u_{n}\circ _{\rho }v_{n};$ then, as $\varepsilon <1/\sigma ,$%
\begin{eqnarray*}
S(z_{n}) &\leq &S(u_{n})\circ _{\sigma }S(v_{n})=S(u_{n})+S(v_{n})(1+\sigma
S(u_{n})) \\
&\leq &\varepsilon (2+\sigma \varepsilon )<3\varepsilon <\lambda _{+}/2.
\end{eqnarray*}%
So%
\[
3\lambda _{+}/4=\lambda _{+}-\varepsilon <S(z_{n})\leq S(u_{n})\circ
_{\sigma }S(v_{n})\leq 3\varepsilon <\lambda _{+}/2,
\]%
a contradiction. So $S$ is continuous at $0$ and so continuous everywhere
(as in Theorem G1):%
\[
-_{\sigma }S(-_{\rho }h)\leq S(x+_{\rho }h)-_{\sigma }S(x)\leq S(h).
\]
The last part follows since $\Sigma \cap (0,\delta ),$ being
Baire/measurable and non-negligible, has the SW property for each $\delta >0$%
. $\square $

\bigskip

\noindent \textbf{Theorem G3.} \textit{Let} $\Sigma \subseteq \lbrack
0,\infty )$ \textit{be locally SW accumulating at }$0$\textit{. Suppose }$S:%
\mathbb{R}\rightarrow \mathbb{R}$\textit{\ is subadditive with }$S(0)=0$
\textit{and }$S|\Sigma $ \textit{is additively bounded above by }$%
G(x):=K_{\kappa }(x)$\textit{, i.e.} $S(\sigma )\leq K_{\kappa }(\sigma )$
\textit{for some }$\kappa $\textit{\ and all }$\sigma \in \Sigma ,$ \textit{%
so that in particular,}%
\[
\lim \sup\nolimits_{\sigma \downarrow 0,\text{ }\sigma \in \Sigma }S(\sigma
)\leq 0.
\]%
\textit{\ Then }$S(x)\leq K_{\kappa }(x)$ \textit{for all }$x>0,$\textit{\ so%
}
\[
\lim \sup_{x\downarrow 0}S(x)\leq 0,
\]%
\textit{and so }$S(0+)=0.$

\textit{In particular, if furthermore there exists a sequence }$%
\{z_{n}\}_{n\in \mathbb{N}}$\textit{\ with }$z_{n}\uparrow 0$\textit{\ and }$%
S(z_{n})\rightarrow 0,$\textit{\ then }$S$\textit{\ is continuous at }$0$%
\textit{\ and so everywhere.}

\bigskip

\noindent \textbf{Proof of Theorem G3.} We are to show that $S(t)\leq
K_{\kappa }(t)$ for all $t.$ We may begin with the simplifying assumption
that $K\equiv 1_{\sigma }=0,$ i.e. that $\kappa =0$, since $S^{\prime
}(t):=S(t)\circ _{\sigma }(K_{\kappa }(t))_{\sigma }^{-1}$ is linearly
bounded above by $1_{\sigma }=0$ on $\Sigma ,$ and $S^{\prime }$ is
subadditive, as $K\ $is additive:%
\[
S^{\prime }(u\circ v):=S(u\circ v)\circ _{\sigma }K_{\kappa }(u\circ _{\rho
}v)_{\sigma }^{-1}\leq S(u)\circ _{\sigma }S(v)\circ _{\sigma }K(u)_{\sigma
}^{-1}\circ _{\sigma }K(v)_{\sigma }^{-1}.
\]%
From now on the proof follows that of [BinO10, Th. 0$^{+}$], mutatis
mutandis (interpreting $+$ as $+_{\rho }$ and $-$ as $-_{\rho }$ as in
Theorem 2). $\square $

\section{Functional inequalities from asymptotic actions: the Goldie argument%
}

We return to the \textit{Karamata asymptotic operator }$K$ acting on $f:%
\mathbb{R}_{+}\rightarrow \mathbb{R}_{+},$ as in ($K$) of \S\ 3, but we now
apply a natural alternative to the limits of \S\ 3 when they cannot be
assumed to exist. This is provided by the $\lim \sup $ operation, which in
the Karamata setting is given by%
\[
K^{\ast }(t)f:=\lim \sup K(t,x)f=\lim \sup \frac{f(xt)}{f(x)}:=K_{f}^{\ast
}(t).
\]%
This leads to an \textit{operator} \textit{domain} defined by%
\[
\mathbb{A}_{f}:=\{u:K_{f}(u):=\lim f(xt)/f(x)\text{ exists and is finite\}.}
\]%
This is a subgroup of $\mathbb{R}_{+}.$ For positive functions $f,$ one has%
\[
\lim \sup \frac{f(xst)}{f(x)}\leq \lim \sup \frac{f(xst)}{f(xt)}\cdot \lim
\sup \frac{f(xt)}{f(x)}:\qquad K^{\ast }(st)f\leq K^{\ast }(s)f\cdot K^{\ast
}(t)f,
\]%
as $K(st,x)\leq K(s,xt)K(t,x).$ Here the limsup yields the multiplicative
\textit{Cauchy functional inequality},%
\begin{equation}
K_{f}^{\ast }(st)\leq K_{f}^{\ast }(s)K_{f}^{\ast }(t),  \tag{$CFI$}
\end{equation}%
as well as a pair of equations restricted to $\mathbb{A}_{f}:$%
\[
\left.
\begin{array}{c}
K_{f}(st)=K_{f}(s)K_{f}(t) \\
K_{f}(t)=K_{f}^{\ast }(t)%
\end{array}%
\right\} \qquad (s,t\in \mathbb{A}_{f}).
\]%
One seeks side-conditions on $f$ and imposes a density condition on $\mathbb{%
A}_{f}$ to deduce that $\mathbb{A}_{f}=\mathbb{R}_{+}.$

For the general asymptotics, with $\varphi \in SE,$\textit{\ }%
\[
K_{h\varphi }^{\ast }(t)f:=\lim \sup \frac{f(x+t\varphi (x))-f(x)}{h(x)},
\]%
there is a corresponding \textit{operator} \textit{domain} defined by%
\[
\mathbb{A}_{hf}:=\{u:K_{hf}(u):=\lim [f(x+t\varphi (x))-f(x)]/h(x)\text{
exists and is finite\},}
\]%
(with $\varphi $ omitted when clear from context). As before, there is also
a functional inequality:%
\[
K_{hf}^{\ast }(t+s\eta (t))\leq K_{hf}^{\ast }(s)K_{h}(t)+K_{hf}^{\ast
}(t),\quad \text{with }K_{h}(t):=\lim h(x\circ _{\varphi }t)/h(x),
\]%
where $K_{h}$ is assumed to exist for all $t$ (as in Prop. 4). The
inequality may be reformulated in Popa-group language as the \textit{%
Beurling-Goldie inequality} satisfied by $K_{f}^{\ast }:\mathbb{G}_{\eta
}\rightarrow \mathbb{G}_{\sigma }:$%
\begin{equation}
K_{hf}^{\ast }(t\circ _{\eta }s)\leq K_{hf}^{\ast }(t)\circ _{\sigma
}K_{hf}^{\ast }(s),\qquad \text{for }\sigma (z)=K_{h}(K_{hf}^{\ast -1}(z)).
\tag{$BGI$}
\end{equation}

However, there is no immediate justification for $\mathbb{A}_{hf}$ being a
subgroup, short of further hypotheses. Either an imposition of good
behaviour of the \textit{limit}, such as local uniformity in $t$, is needed,
thus \textit{narrowing} the domain, or a presumption of topologically good
character of the \textit{domain} itself, such as requiring $\mathbb{A}_{hf}$
to contain a non-meagre subset. The latter may draw on additional axioms of
set theory, for which see [BinO12]. For an extensive study of the uniformity
assumptions, see [BinO7].

Henceforth we take for granted a domain $\mathbb{A}$ that is a dense
subgroup of an appropriate Popa group $\mathbb{G}$, and a side-condition of
right-sided continuity at $0=1_{\mathbb{G}}$ imposed on $K_{hf}^{\ast }$ (so
on $\mathbb{R}_{+}$).

Above we had the Beurling-Goldie equation ($BG$). Below, we restrict one or
both of the arguments $u$ and $v$ to $\mathbb{A}$, obtaining the `singly
conditioned' and `doubly conditioned' Beurling-Goldie equations ($BG_{%
\mathbb{A}}$) and ($BG_{\mathbb{AA}}$). For the origins of the \textit{%
Goldie argument}, see the Remark after Theorem 4 below.

We begin with an auxiliary result. (In the equation below $g(0)K(0)=0$, so
to avoid trivial (constant) solutions w.l.o.g. we assume both here and later
that $g(0)=1.)$

\bigskip

\noindent \textbf{Proposition 6 }([BojK, (2.2)], [BinGT, Lemma 3.2.1]; cf.
[AczG]). \textit{Take }$\eta \in GS$ \textit{and }$g$ \textit{with }$g(0)=1.$
\textit{If }$K\not\equiv 0$ \textit{satisfies}%
\begin{equation}
K(u\circ _{\eta }v)=g(v)K(u)+K(v)\text{\qquad }(u,v\in \mathbb{A}),
\tag{$BG_{\QTR{Bbb}{AA}}$}
\end{equation}%
\textit{\ \ with }$\mathbb{A}$\ \textit{a dense subgroup of }$\mathbb{G}%
_{\eta }$\textit{,} \textit{then:\newline
}\noindent (i)\textit{\ the following is a subgroup of }$\mathbb{G}_{\eta }$
\textit{on which }$K$\textit{\ is additive:}%
\[
\mathbb{A}_{g}:=\{u\in \mathbb{A}:\mathit{\ }g(u)=1\};
\]%
\textit{\ }\noindent (ii) \textit{if }$\mathbb{A}_{g}\neq \mathbb{A}$
\textit{and }$K\not\equiv 0$$,$ \textit{there is a constant }$\kappa \neq 0$
\textit{with}
\begin{equation}
K(t)\equiv \kappa (g(t)-1)\text{\qquad }(t\in \mathbb{A}),  \tag{*}
\end{equation}%
\textit{\ and }$g$\textit{\ satisfies }%
\begin{equation}
g(u\circ _{\eta }v)=g(v)g(u)\text{\qquad }(u,v\in \mathbb{A}).  \tag{$CJ$}
\end{equation}%
\noindent (iii) \textit{So for }$\mathbb{A=G}_{\eta }^{+}$\textit{\ with }$%
\eta =\eta _{\rho }$ \textit{and }$g$\textit{\ locally bounded at }$0$%
\textit{\ with }$g\neq 1$\textit{\ except at }$0:$\textit{\ }%
\[
g(x)\equiv (1+\rho t)^{\gamma },
\]%
\textit{for some constant }$\gamma \neq 0,$ \textit{and so }$K(t)\equiv
cK_{\gamma }(t)$ \textit{for some constant }$c$, \textit{where}%
\[
K_{\gamma }(t):=[(1+\rho t)^{\gamma }-1].
\]

\bigskip

\noindent \textit{Proof.} This is proved exactly as in [BinO6, Th. 1] with $%
\circ _{\eta _{\rho }}$or $+_{\rho }$ replacing $+$. One uses the \textit{%
Cauchy nucleus }of $K$ [Kuc, Lemma 18.5.1]. $\square $

\bigskip

\noindent \textbf{Example }in the case $\rho =1$\textbf{.} Below, put $x=u+1$
and $k(t)=g(t-1)$ :
\begin{eqnarray*}
g((u+1)(v+1)-1) &=&g(u+v+uv)=g(u)g(v):\qquad g(xy-1)=g(x-1)g(y-1); \\
k(xy) &=&k(x)k(y):\qquad g(t)=k(t+1)=(1+t)^{\gamma }; \\
K(t) &=&\kappa (g(t)-1)=\kappa \lbrack (1+t)^{\gamma }-1].
\end{eqnarray*}

\bigskip

\noindent \textbf{Theorem 4 (Generalized Goldie Theorem, }cf. [BinO6, Th. 3]%
\textbf{).} \textit{If for }$\eta \in GS$ \textit{and }$\mathbb{A}$\textit{\
a dense subgroup of }$\mathbb{G}_{\eta }$\textit{\ ,}

\noindent (i)\textit{\ }$F^{\ast }:$\textit{\ }$\mathbb{R\rightarrow R}$
\textit{is positive and subadditive with }$F^{\ast }(0+)=0;$

\noindent (ii) $F^{\ast }$ \textit{satisfies the singly-conditioned
Beurling-Goldie equation}%
\begin{equation}
F^{\ast }(u\circ _{\eta }v)=g(v)K(u)+F^{\ast }(v)\text{\qquad }(u\in \mathbb{%
A})(v\in \mathbb{R}_{+})  \tag{$BG_{\QTR{Bbb}{A}}$}
\end{equation}%
\textit{for some non-zero }$K$\textit{\ satisfying }$(BG_{\mathbb{A}})$%
\textit{\ with }$g$ \textit{continuous on }$\mathbb{R}$\textit{\ and }$%
\mathbb{A}_{g}=\{0\}$\textit{\ (i.e. }$g(0)=1$ \textit{but otherwise }$%
g(v)\neq 1);$

\noindent (iii)\textit{\ }$F^{\ast }$ \textit{extends }$K$ \textit{on }$%
\mathbb{A}$:\textit{\ }%
\[
F^{\ast }(x)=K(x)\text{\qquad }(x\in \mathbb{A}),
\]%
\textit{so that in particular }$F^{\ast }$\textit{\ satisfies }$(BG_{\mathbb{%
A}}),$ \textit{and indeed}%
\[
F^{\ast }(u\circ _{\eta }v)=g(v)F^{\ast }(u)+F^{\ast }(v)\text{\qquad }(u\in
\mathbb{A})(v\in \mathbb{G}_{\eta }^{+}):
\]%
-- \textit{then for some }$c>0,\gamma \geq 0$%
\[
g(x)\equiv c(1+\rho t)^{-\gamma }\text{ and }F^{\ast }(x)\equiv cK_{-\gamma
}(x)=c[(1+\rho t)^{-\gamma }-1]/\rho \text{\qquad }(x\in \mathbb{R}_{+}).
\]

\bigskip

\noindent \textit{Proof. }We write $\circ $ for $\circ _{\eta }$, and $%
\mathbb{G}$ for $\mathbb{G}_{\eta }.$ Put%
\[
G(x):=\int_{0}^{x}g(t)\text{ }\mathrm{d}t/\eta (t):\text{\qquad }G^{\prime
}(x)=g(x)/\eta (x).
\]%
By continuity of $g$ and Th. 1, $K$ is continuous on $\mathbb{A}$, so $%
K(u+)=K(u)$ for all $u\in \mathbb{A}$, and so in particular $K(0+)=0,$ which
is also implied by (i) above. Also note that $F^{\ast }$ is right-continuous
(and $F^{\ast }(u+)=K(u)$) on $\mathbb{A}$, and on $\mathbb{G}$ satisfies%
\[
\lim \sup_{v\downarrow 0}F^{\ast }(u\circ v)\leq g(0)F^{\ast }(u)+F^{\ast
}(0+)=F^{\ast }(u).
\]

We write $\delta ^{n\circ }$ for the $n$\textit{-fold} \textit{product} in $%
\mathbb{G}$ (inductively defined so that $\delta ^{0\circ }=1_{\mathbb{G}}=0$
and $\delta ^{n\circ }=\delta ^{(n-1)\circ }\circ \delta ).$

Now we mimick the Goldie proof of [BinGT, \S 3.2.1] (extending [BinO6, Th.
3] to the current Popa context). For any $u,u_{0}$ with $u_{0}\in \mathbb{A}$
and $u_{0}>0$, define $i=i(\delta )\in \mathbb{Z}$ for $\delta >0$ so that $%
\delta ^{(i-1)\circ }\leq u<\delta ^{i\circ },$ and likewise for $u_{0}$
define $j=i_{0}(\delta ).$ As $\mathbb{A}_{g}=\{0\},$ put $%
c_{0}:=K(u_{0})/[g(u_{0})-1].$ For $m\in \mathbb{N}$%
\[
F^{\ast }(\delta ^{m\circ })-F^{\ast }(\delta ^{(m-1)\circ })=g(\delta
^{(m-1)\circ })K(\delta ),
\]%
as $\delta ^{m\circ }\in \mathbb{A}$, so that on summing%
\begin{equation}
F^{\ast }(\delta ^{i\circ })=K(\delta )\sum_{m=1}^{i}g(\delta ^{(m-1)\circ
}),  \tag{**}
\end{equation}%
as $F^{\ast }(0)=0.$ Note that%
\[
\Delta _{m}:=\delta ^{m\circ }-\delta ^{(m-1)\circ }=\delta \eta (\delta
^{(m-1)\circ }),
\]%
so that $\Delta _{m}\rightarrow 0$ as $\delta \rightarrow 0.$ But as $%
G^{\prime }(t)=g(t)/\eta (t),$%
\begin{equation}
\sum_{m=1}^{i}g(\delta ^{m\circ })\Delta _{m}=\sum_{m=1}^{i}G^{\prime
}(\delta ^{(m)\circ })\eta (\delta ^{m\circ })\delta \rightarrow
\int_{0}^{u}G^{\prime }(x)\text{ }\mathrm{d}x  \tag{$RI$}
\end{equation}%
(for `Riemann Integral'). Without loss of generality $G(u_{0})\neq 0.$
(Indeed, otherwise $g=0$ on $\mathbb{A}\cap \mathbb{R}_{+}$ and so on $%
\mathbb{R}_{+}$, so that $F^{\ast }(u+)=0$ on $\mathbb{A}\cap \mathbb{R}_{+}$%
; this together with $F^{\ast }(u+v)=F^{\ast }(v)$ contradicts positivity of
$F^{\ast }$ on $\mathbb{R}_{+}$.) Taking limits as $\delta \rightarrow 0$
through positive $\delta \in \mathbb{A}$ with $K(\delta )\neq 0$ (see below
for $K(\delta )=0$), we then have, as $G(u_{0})\neq 0,$%
\begin{eqnarray*}
\frac{F^{\ast }(\delta ^{i\circ })}{F^{\ast }(\delta ^{j\circ })} &=&\frac{%
K(\delta )}{K(\delta )}\frac{\sum_{m=1}^{i}g(\delta ^{m\circ })}{%
\sum_{m=1}^{j}g(\delta ^{m\circ })}=\frac{\sum_{m=1}^{i}G^{\prime }(\delta
^{(m)\circ })\eta (\delta ^{m\circ })\delta }{\sum_{m=1}^{i_{0}}G^{\prime
}(\delta ^{(m)\circ })\eta (\delta ^{m\circ })\delta }=\frac{%
\sum_{m=1}^{i}G^{\prime }(\delta ^{(m)\circ })\Delta _{m}}{%
\sum_{m=1}^{i_{0}}G(\delta ^{(m)\circ })\Delta _{m}} \\
&\rightarrow &\frac{\int_{0}^{u}G^{\prime }(x)\text{ }\mathrm{d}x}{%
\int_{0}^{u_{0}}G^{\prime }(x)\text{ }\mathrm{d}x}=\frac{G(u)}{G(u_{0})}.
\end{eqnarray*}%
Here by right-continuity at $u_{0}$
\[
\lim F^{\ast }(\delta ^{i\circ })=F^{\ast
}(u_{0})=K(u_{0})=c_{0}[g(u_{0})-1]>0.
\]%
So%
\[
F^{\ast }(\delta ^{i\circ })\rightarrow G(u)\cdot F^{\ast }(u_{0})/G(u_{0}).
\]%
Put $c_{1}:=c_{0}[g(u_{0})-1]/G(u_{0}).$ As before, as $u_{0}\in \mathbb{A},$%
\begin{eqnarray*}
F^{\ast }(u) &\geq &\lim \sup F^{\ast }(\delta ^{i\circ })=G(u)\cdot F^{\ast
}(u_{0})/G(u_{0}) \\
&=&G(u)K(u_{0})/G(u_{0})=G(u)c_{0}[g(u_{0})-1]/G(u_{0})=c_{1}G(u).
\end{eqnarray*}%
Now specialize to $u\in \mathbb{A},$ on which, by above, $F^{\ast }$ is
right-continuous. Letting $\delta ^{i\circ }\in \mathbb{A}$ decrease to $u,$
the inequality above becomes an equation:%
\[
K(u)=F^{\ast }(u)=c_{1}G(u)\text{\qquad }(u\in \mathbb{A}).
\]%
(This remains valid with $c_{1}=0$ if $K(\delta )=0$ for $\delta \in \mathbb{%
A}\cap I$ for some interval $I=(0,\varepsilon ),$ since then $F^{\ast }(u)=0$
by right-continuity on $\mathbb{A}$, as $F^{\ast }(\delta ^{i\circ })=0$ for
$\delta \in \mathbb{A}\cap I,$ by (**).)

We extend the domain of this equation from $\mathbb{A}$ to the whole of $%
\mathbb{R}$, using a key idea due to Goldie (see the Remark below).

For an \textit{arbitrary} $u\in \mathbb{R},$ take $v\in \mathbb{A}$ with $%
z:=u-v>0,$ i.e. with $v<u.$ Then%
\begin{eqnarray*}
F^{\ast }(u) &=&F^{\ast }(v+z)=K(v)g(z)+F^{\ast }(z)\qquad \text{(by (ii),
as }v\in \mathbb{A}\text{ and }z\in \mathbb{R}_{+}\text{)} \\
&=&c_{1}G(v)g(z)+F^{\ast }(z) \\
&\rightarrow &c_{1}G(u)g(0)+0=c_{1}G(u)\qquad \text{(as }z\downarrow 0\text{%
),}
\end{eqnarray*}%
by continuity of $g$ and $G,$ and $F^{\ast }(0+)=0.$ So%
\[
F^{\ast }(u)=c_{1}G(u)\text{\qquad }(u\in \mathbb{G}).
\]%
Thus by (*) of Prop. 6, for some $\kappa $%
\[
c_{1}G(u)=F^{\ast }(u)=K(u)=\kappa \lbrack g(u)-1]\text{\qquad }(u\in
\mathbb{A})\text{.}
\]%
So, by density and continuity on $\mathbb{G}$ of $g,$
\[
\kappa \lbrack g(u)-1]=c_{1}G(u)\text{\qquad }(u\in \mathbb{R}_{+}).
\]%
Thus $g$ is indeed differentiable; differentiation now yields%
\[
\kappa g^{\prime }(u)=c_{1}g(u)/\eta (u):\qquad g^{\prime
}(u)/g(u)=(c_{1}/\kappa \eta (u))\text{\qquad }(u\in \mathbb{R}_{+}),
\]%
as $\kappa \neq 0$ (otherwise $K(u)\equiv 0,$ contrary to assumption). So,
as $g(0)=1,$ with $\gamma :=-c_{1}/\kappa \rho $%
\[
\log g(u)=-\frac{c_{1}}{\kappa }\int_{0}^{u}\frac{\mathrm{d}t}{1+\rho t}%
=-\gamma \log (1+\rho u):\text{\qquad }g(u)=(1+\rho t)^{-\gamma }.
\]%
So%
\[
G(u)=\int_{0}^{u}g(t)\frac{\mathrm{d}t}{\eta (t)}=\int_{0}^{u}(1+\rho
t)^{-\gamma -1}\mathrm{d}t=[(1+\rho u)^{-\gamma }-1]/\rho .
\]%
So%
\[
\text{ }G(u)=cK_{\gamma }(u):\text{\quad }F^{\ast
}(u)=c_{1}G(u)=c_{1}[(1+\rho u)^{-\gamma }-1]/\rho \text{\quad }(u\in
\mathbb{R}).
\]%
As $(1+\rho u)^{-\gamma }$ is subadditive on $\mathbb{R}_{+}$ iff $\gamma
\geq 0$ (cf. before Th. 1), $c_{1}>0.$ $\square $

\bigskip

\noindent \textbf{Remark. }Above, we have disaggregated the Goldie proof
given in [BinGT, \S 3.2.1] into three steps. Firstly, we use the integral $G$
of the unknown auxiliary function $g$ (as in [BinO6, Th. 3], albeit here as
a Haar integral), where Goldie assumed $g$ explicitly to be the exponential
function $e^{\gamma t}.$ For Goldie this permits an explicit formula for the
corresponding sums (for us the Riemann sums lead to a simple differential
equation, which we can solve for $g$, giving $G$). Secondly, we have
partitioned the range of integration by use of a \textit{Beck sequence}
[Bec, Lemma 1.64] (iterating $\circ \delta $). Finally, the extension of the
relation between $F^{\ast }$ and $G$ from $\mathbb{A}$ to $\mathbb{R}_{+}$
makes explicit a remarkable achievement, due to Goldie (and left implicit in
[BinGT, \S\ 3.2.1]): establishment of left-sided continuity from the assumed
right-sided continuity $F^{\ast }(0+)=0.$ This overlooked feature was first
made explicit in [BinO10] as Theorem 0 there (cf. Th. G1 above), yielding
new results, and again put to further extensive use in [BinO11].

\bigskip

Armed with the results here we are now able to freely lift results from
[BinO10] concerning when the solution $K_{hf}^{\ast }:\mathbb{G}_{\eta
}\rightarrow \mathbb{G}_{\sigma }$ of ($BGI$) in fact solves ($BG$) and so
takes the form $K_{\kappa }(u)$\textit{\ for some} $\kappa \in \mathbb{R}$.
We recall that in the interests of simplicity we assume that the domain of
the asymptotic operator is a subgroup, leaving the reader to refer for
results which guarantee this to [BinO7]. Below, we use \textit{linear} to
mean continuous and additive.

\bigskip

\noindent \textbf{Theorem 5 (Quantifier-Weakening Theorem, }cf. [BinO10, Th.
6], [BinO7, Th. 6]\textbf{).} \textit{With }$K_{hf}^{\ast }$ \textit{and }$%
\mathbb{A}_{hf}$\textit{\ as above, suppose that}\newline
\noindent (i)\textit{\ }$\mathbb{A}_{hf}$ \textit{is a dense subgroup of }$%
\mathbb{G}_{\eta }$\textit{;}\newline
\noindent (ii)\textit{\ }$K_{hf}^{\ast }$\textit{\ satisfies the one-sided
Heiberg-Seneta boundedness condition }%
\begin{equation}
\lim \sup\nolimits_{u\downarrow 0}K_{hf}^{\ast }(u)\leq 0  \tag{$HS$}
\end{equation}%
-- \textit{then }$\mathbb{A}_{hf}=\mathbb{G}_{\eta }$\textit{\ and }$%
K_{hf}^{\ast }$\textit{\ is linear (continuous and additive):}%
\[
K_{hf}^{\ast }(u)=\lim_{x\rightarrow \infty }[f(u+x)-f(x)]/h(x)=K_{\kappa
}(u)
\]%
\textit{\ for some} $\kappa \in \mathbb{R}$\textit{, and all }$u\in \mathbb{G%
}_{\eta }.$

\bigskip

\noindent \textbf{Proof of Theorem 5. }As we assume here that $\mathbb{A}%
_{hf}$ is a subgroup, referring to results in [BinO10, Props 3 and 6], $%
K_{f}^{\ast }$ is a finite, subadditive, right-continuous extension of $G$.
So $G$ is continuous on $\mathbb{A}_{f},$ and so $G(\sigma )=K_{\kappa
}(\sigma ),$ for all $\sigma \in \mathbb{A}_{f}.$ As $\mathbb{A}_{f}$ is
dense, by [BinO10, Prop. 7], $K_{f}^{\ast }(u)=K_{\kappa }(u)$ for all $u.$
By [BinO10, Prop. 1], $\mathbb{A}_{f}=\mathbb{G}$ and $K_{f}^{\ast
}(u)=G(u). $ $\square $

\bigskip

We turn now to\textit{\ thinnings} of the condition ($HS$) of Theorem 5. For
this we need some definitions from [BinO10].

\bigskip

\noindent \textbf{Definitions.} 1. Say that $\Sigma $ is \textit{locally
Steinhaus-Weil (SW)}, or has the \textit{SW property locally}, if for $%
x,y\in \Sigma $ and, for all $\delta >0$ sufficiently small, the sets%
\[
\Sigma _{z}^{\delta }:=\Sigma \cap B_{\delta }(z),
\]%
for $z=x,y,$ have the \textit{interior-point property,} that $\Sigma
_{x}^{\delta }\pm \Sigma _{y}^{\delta }$ has $x\pm y$ in its interior. (Here
$B_{\delta }(x)$ is the open ball about $x$ of radius $\delta .)$ See
[BinO3] for conditions under which this property is implied by the
interior-point property of the sets $\Sigma _{x}^{\delta }-\Sigma
_{x}^{\delta }$ (cf. [BarFN]); see also the rich list of examples below,
which are used in [BinO8,10,11,13,14], [MilMO].

\noindent 2. Say that $\Sigma \subseteq \mathbb{R}$ is \textit{shift-compact
}if for each \textit{null sequence} $\{z_{n}\}$ (i.e. with $z_{n}\rightarrow
0$) there are $t\in \Sigma $ and an infinite $\mathbb{M\subseteq N}$ such
that%
\[
\{t+z_{m}:m\in \mathbb{M}\}\subseteq \Sigma .
\]%
See [BinO4], and for the group-action aspects, [MilO].

\bigskip

\noindent \textbf{Examples} \textbf{of families of locally Steinhaus-Weil
sets} (see e.g. [BinO13]).

The sets listed below are typically, though not always, members of a
topology on an underlying set.

\noindent (o) $\Sigma $ a usual (Euclidean) open set in $\mathbb{R}$ (and in
$\mathbb{R}^{n})$ -- this is the `trivial' example;

\noindent (i) $\Sigma $ density-open subset of $\mathbb{R}$ (similarly in $%
\mathbb{R}^{n})$ (by Steinhaus's Theorem -- see e.g. [BinGT, Th. 1.1.1],
[BinO13], [Oxt, Ch. 8]);

\noindent (ii) $\Sigma $ locally non-meagre at all points $x\in \Sigma $ (by
the Piccard-Pettis Theorem -- as in [BinGT, Th. 1.1.2], [BinO13], [Oxt, Ch.
8] -- such sets can be `thinned out', i.e. extracted as subsets of a
second-category set, using separability or by reference to the Banach
Category Theorem [Oxt, Ch.16]);

\noindent (iii) $\Sigma $ the Cantor `excluded middle-thirds' subset of $%
[0,1]$ (since $\Sigma +\Sigma =[0,2]);$

\noindent (iv) $\Sigma $ universally measurable and open in the \textit{ideal%
} topology ([LukMZ], [BinO9]) generated by omitting Haar null sets (by the
Christensen-Solecki Interior-points Theorem of [Sol]);

\noindent (v) $\Sigma $ a Borel subset of a Polish abelian group and and
open in the ideal topology generated by omitting \textit{Haar meagre} sets
in the sense of Darji [Dar] (by Jab\l o\'{n}ska's generalization of the
Piccard Theorem, [Jab1, Th. 2], cf. [Jab3], and since the Haar-meagre sets
form a $\sigma $-ideal [Dar, Th. 2.9]); for details see [BinO13].

If $\Sigma $ is \textit{Baire} (has the Baire property) and is locally
non-meagre, then it is co-meagre (since its quasi interior is everywhere
dense).

\noindent \textbf{Caveat. }1. Care is needed in identifying locally SW sets:
Mato\u{u}skov\'{a} and Zelen\'{y} [MatZ] show that in any non-locally
compact abelian Polish group there are closed non-Haar null sets $A,B$ such
that $A+B$ has empty interior. Recently, Jab\l o\'{n}ska [Jab4] has shown
that likewise in any non-locally compact abelian Polish group there are
closed non-Haar meager sets $A,B$ such that $A+B$ has empty interior.

\noindent 2. For an example on $\mathbb{R}$ of a compact subset $S$ such
that $S-S$ does contains an interval, but $S+S$ has measure zero and so does
not, see [CrnGH].

\noindent 3. Here we are concerned with subsets $\Sigma \subseteq \mathbb{R}$
where such `anomalies' are assumed not to occur.

\bigskip

We can now state some thinned variants of Th. 6.

\bigskip

\noindent \textbf{Theorem 6 (Thinned Quantifier Weakening Theorem}; [BinO10,
Th. 1$^{\prime }$], cf. [BinO7, \S 6 Th. 5]). \textit{Theorem 5 above holds
with condition }(ii)\textit{\ replaced by any one of the following:}\newline
\noindent (ii-a)\textit{\ }$K_{hf}^{\ast }$\textit{\ satisfies the
Heiberg-Seneta boundedness condition thinned out to a symmetric set }$\Sigma
$\textit{\ that is locally SW, i.e.}%
\[
\lim \sup\nolimits_{u\rightarrow 0,\text{ }u\in \Sigma }K_{hf}^{\ast
}(u)\leq 0;
\]%
\noindent (ii-b)\textit{\ }$K_{hf}^{\ast }$\textit{\ is linearly bounded
above on a locally SW subset }$\Sigma \subseteq $\textit{\ }$\mathbb{R}%
_{+}=(0,\infty )$\textit{\ accumulating at }$0,$\textit{\ so that in
particular}%
\[
\lim \sup\nolimits_{u\downarrow 0,\text{ }u\in \Sigma }K_{hf}^{\ast }(u)\leq
0;
\]%
\textit{\ }\noindent (ii-c)\textit{\ }$K_{hf}^{\ast }$\textit{\ is bounded
above on a locally SW subset }$\Sigma \subseteq $\textit{\ }$\mathbb{A}_{+}$%
\textit{\ accumulating at }$0$\textit{, that is, the following }$\lim \sup $
\textit{is finite:}

\begin{equation}
\lim \sup\nolimits_{u\downarrow 0,\text{ }u\in \Sigma }K_{hf}^{\ast
}(u)<\infty ;  \tag{$SW$-$HS(K_{hf}^{\ast })$}
\end{equation}%
\noindent (ii-d) $S$\textit{\ is bounded on a subset }$\Sigma \subseteq $%
\textit{\ }$\mathbb{A}$\ \textit{that is shift-compact (e.g. on a set that
is locally SW, and so on an open set) and\ }%
\[
\mathbb{A}=\mathbb{A}_{hf}:=\{u:K_{hf}(u):=\lim_{x\rightarrow \infty
}[f(u+x)-f(x)]/h(x)\text{ exists and is finite}\}.
\]

\noindent \textbf{Proof. }This follows from the Popa variant of [BinO10,
Theorem 1$^{\prime }$], the proof of which follows from Theorems G2 and G3
of \S 4 above in place of [BinO10, Theorems 0$^{\prime }$ and 0]. $\square $

\bigskip

The classical \textit{Quantifier Weakening Theorems} of regular variation
([BinGT,\ \S 1.4.3 and \S 3.2.5]) are re-stated below as Theorems K and
BKdH. There, one needs as side-condition the Heiberg-Seneta condition $HS$
restated multiplicatively here as ($\lim \sup $) (or a thinned version of
it, as in Theorem 6). Recall from above the $^{\ast }$ notation (as in $%
g^{\ast })$ signifying that limsup replaces $\lim .$

\bigskip

\noindent \textbf{Theorem K} (cf. [BinGT, Th. 1.4.3]). \textit{Suppose that}%
\begin{equation}
\lim \sup\nolimits_{\lambda \downarrow 1}K_{f}^{\ast }(\lambda )\leq 1.
\tag{$\lim \sup $-$f$}
\end{equation}%
\textit{Then the following are equivalent:}\newline
\textit{\noindent }(i)\textit{\ there exists }$\rho \in \mathbb{R}$\textit{\
such that }%
\[
f(\lambda x)/f(x)\rightarrow \lambda ^{\rho }\qquad (x\rightarrow \infty
)(\forall \lambda >0);
\]%
\textit{\noindent }(ii)\textit{\ }$g(\lambda )=\lim_{x\rightarrow \infty
}f(\lambda x)/f(x)$\textit{\ exists, finite for all }$\lambda $\textit{\ in
a non-negligible set;}\newline
\textit{\noindent }(iii)\textit{\ }$g(\lambda )$\textit{\ exists, finite,
for all }$\lambda $\textit{\ in a dense subset of }$(0,\infty );$\newline
\textit{\noindent }(iv)\textit{\ }$g(\lambda )$\textit{\ exists, finite for }%
$\lambda =\lambda _{1},\lambda _{2}$\textit{\ with }$(\log \lambda
_{1})/\log \lambda _{2}$ \textit{irrational.}

\bigskip

Theorem K\ is an immediate corollary of Theorem 5, as (limsup) iff ($%
HS(K_{f}^{\ast })$). The final assertion follows from Kronecker's theorem
[HarW, Ch. 23].

\bigskip

\noindent \textbf{Theorem BKdH} (cf. [BinGT, Th. 3.2.5]). \textit{For }$h$
\textit{with }%
\[
\lim_{x\rightarrow \infty }h(\lambda x)/h(x)=\lambda ^{\rho }\qquad (\lambda
>0),
\]%
\textit{\ and}%
\begin{equation}
\lim \sup\nolimits_{\lambda \downarrow 1}K_{hf}^{\ast }(\lambda )\leq 0,
\tag{$\lim \sup $-$hf$}
\end{equation}%
\textit{the following are equivalent:}\newline
\textit{\noindent }(i)\textit{\ }$K_{hf}(\lambda ):=\lim_{x\rightarrow
\infty }$[$f(\lambda x)-f(x)]/h(x)$\textit{\ exists, finite for all }$%
\lambda >0,$ \textit{and }$K_{hf}(\lambda )=c\eta _{\rho }^{-1}(\lambda
^{\rho })$\textit{\ for some }$c$\textit{\ and all }$\lambda $\textit{\ on a
non-negligible set;}\newline
\textit{\noindent }(ii)\textit{\ }$K_{hf}(\lambda )$ \textit{exists, finite
for all }$\lambda $\textit{\ in a non-negligible set;}\newline
\textit{\noindent }(iii)\textit{\ }$K_{hf}(\lambda )$\textit{\ exists,
finite, for all }$\lambda $\textit{\ in a dense subset of }$(0,\infty );$%
\newline
\textit{\noindent }(iv)\textit{\ }$K_{hf}(\lambda )$\textit{\ exists, finite
for }$\lambda =\lambda _{1},\lambda _{2}$\textit{\ with }$(\log \lambda
_{1})/\log \lambda _{2}$ \textit{irrational.}

\bigskip

Theorem BKdH is an immediate corollary of Theorem 4. As before the final
assertion follows from Kronecker's theorem.

The motivation for this paper was the treatment of Theorems K and BKdH above
via Popa groups in [BinO7, \S 7] (specifically ($GFE$) and ($GS$) there and
their equivalence), using the extra power of the extra generality here to
provide a unified and simplified treatment.

\bigskip

\section{Concluding Remarks}

\noindent \textit{Beurling's Tauberian theorem. }To extend the Wiener
Tauberian Theorem (Theorem W, say) Beurling introduced (in unpublished
lectures of 1957) his Tauberian Theorem (below), extending Theorem W from
convolutions to `convolution-like' operations. We need the \textit{Beurling
convolution}:%
\[
F\ast _{\varphi }H(x):=\int F\Bigl(\frac{x-u}{\varphi (x)}\Bigr)H(u)\frac{%
\mathrm{d}u}{\varphi (x)}=\int F(-t)H(x\circ _{\varphi }t)\text{ }\mathrm{d}%
t.
\]%
This is an asymptotic version, involving the function $\eta _{x}(.)$ of \S 3:%
\[
\eta _{x}(t):=\varphi (x\circ _{\varphi }t)/\varphi (x),
\]%
of an ordinary convolution (below).

\bigskip

\noindent \textbf{Theorem B (Beurling's Tauberian theorem)}.\textit{\ For }$%
K\in L_{1}(\mathbb{R})$\textit{\ with Fourier transform }$\hat{K}$\textit{\
non-zero on }$\mathbb{R}$\textit{, and }$\varphi $ \textit{Beurling slowly
varying, that is}%
\begin{equation}
\eta _{x}(t)\rightarrow 1,\qquad (x\rightarrow \infty )\qquad (t\geq 0):
\tag{$BSV$}
\end{equation}

\textit{if }$H$\textit{\ is bounded, and}%
\[
K\ast _{\varphi }H(x)\rightarrow c\int K(y)\mathrm{d}y,
\]%
\textit{\ then for all }$F\in L_{1}(\mathbb{R})$
\[
F\ast _{\varphi }H(x)\rightarrow c\int F(y)\mathrm{d}y\qquad (x\rightarrow
\infty ).
\]%
This reduces to Theorem W on replacing $\varphi $ by 1. For an elegant
proof, see [Kor, IV.11].

In Theorem W, the argument in the integral above (with $\varphi =1$) is $x-u$%
, and so is a convolution (written additively, or $x/u$ multiplicatively).
In Theorem B, the integral is merely `convolution-like'. Beurling was able
to use his form of slow variation, ($BSV$), to reduce easily to convolution
form, and so to Theorem W. His motivation was the Tauberian theorem for the
Borel summability method, important in summability theory, complex analysis
and probability [Kor, VI]. For applications in probability, see e.g.
[Bin1,3].

Beurling convolution is an \textit{asymptotic convolution}: to within a
factor $\eta _{x}(t)\rightarrow 1$, it is the proper convolution%
\[
(f\ast _{\varphi }g)(x):=\int_{G_{\rho }}f(-t/\eta _{x}(t))g(x\circ
_{\varphi }t)\text{ }\mathrm{d}\eta _{G_{\rho }}(t)\qquad (x\in G_{\rho }).
\]%
For, given $x$ and $t,$ solving for $s$ the equation%
\[
x=(x\circ _{\varphi }t)\circ _{\varphi }s=x+t\varphi (x)+s\varphi
(x+t\varphi (x))
\]%
yields%
\[
s=-t\varphi (x)/\varphi (x+t\varphi (x))=-t/\eta _{x}(t)
\]%
as the `inverse of $t$' (relative to the binary operation $\circ _{\varphi }$
acting on the \textit{set} $G_{\rho }$).

For $\varphi \in SE,$ the corresponding asymptotic convolution is
\[
(f\ast _{\varphi }g)(x):=\int f(-t/\eta _{\rho }(t))g(x\circ _{\varphi }t)%
\text{ }\mathrm{d}\eta _{G_{\rho }}(t).
\]%
For $\varphi (x):=\eta _{\rho }(x)\in GS,$ since%
\[
\eta _{x}(t):=\frac{\varphi (x+t\varphi (x))}{\varphi (x)}=\frac{\eta _{\rho
}(x+t\eta _{\rho }(x))}{\eta _{\rho }(x)}=\eta _{\rho }(t),
\]%
$(f\ast _{\varphi }g)(x)$ becomes
\[
(f\ast _{\eta _{\rho }}g)(x):=\int f(-t/\eta _{\rho }(t))g(x\circ _{\eta
_{\rho }}t)\text{ }\mathrm{d}\eta _{G_{\rho }}(t)=\int f(-_{\rho
}t)g(x+_{\rho }t)\text{ }\mathrm{d}\eta _{G_{\rho }}(t),
\]%
with the notation of Theorem 2. So in this case the asymptotic convolution
becomes ordinary convolution for the Popa group $(G_{\rho },\circ _{\rho }).$

\bigskip

\textbf{Postscript.}

The whole area of regular variation stems from the pioneering work of Jovan
Karamata (1902-1967) in 1930. The present paper stems from his work with
Ranko Bojanic (1925-2017) of 1963 [BojK]. The first author offers here a
reminiscence of his first meeting with Ranko Bojanic (in 1988, over dinner,
at a conference at Ohio State University, Columbus OH). He asked Professor
Bojanic why he and Karamata had stopped their work on regular variation in
1963. He replied unhesitatingly `Because we didn't know what it was good
for'. Analysts in general, and probabilists in particular, do now know what
it is good for. Our aim here has been to demonstrate the power, and ongoing
influence, of their work, with the benefit of 55 years worth of hindsight.

\bigskip

\textbf{References}

\noindent \lbrack AczD] J. Acz\'{e}l and J. Dhombres, \textsl{Functional
equations in several variables.} Encycl. Math. App.\textbf{\ 31}, Cambridge
University Press, 1989.\newline
\noindent \lbrack AczG] J. Acz\'{e}l and S. Go\l \k{a}b S, Remarks on
one-parameter subsemigroups of the affine group and their homo- and
isomorphisms. \textsl{Aequ. Math.} \textbf{4 }(1970), 1--10.\newline
\noindent \lbrack ArhT] A. Arhangel'skii and M. Tkachenko, \textsl{%
Topological groups and related structures}. World Scientific, 2008.\newline
\noindent \lbrack BarFN] A. Bartoszewicz, M. Filipczak and T. Natkaniec, On
Sm\'{\i}tal properties. \textsl{Topology Appl.} \textbf{158} (2011),
2066--2075.\newline
\noindent \lbrack Bec] A. Beck, \textsl{Continuous flows on the plane},
Grundl. math. Wiss. \textbf{201}, Springer, 1974.

\noindent \lbrack Bil] P. Billingsley, \textsl{Ergodic theory and information%
}. Wiley, 1965.\newline
\noindent \lbrack BinG] N. H. Bingham and C. M. Goldie. On one-sided
Tauberian conditions, \textsl{Analysis} \textbf{3} (1983), 159-188.\newline
\noindent \lbrack BinGT] N. H. Bingham, C. M. Goldie and J. L. Teugels,
\textsl{Regular variation}, 2nd ed., Cambridge University Press, 1989 (1st
ed. 1987).\newline
\noindent \lbrack BinO1] N. H. Bingham and A. J. Ostaszewski, Beyond
Lebesgue and Baire II: bitopology and measure-category duality. \textsl{%
Colloq. Math.} \textbf{121} (2010), no. 2, 225--238.\newline
\noindent \lbrack BinO2] N. H. Bingham and A. J. Ostaszewski, Normed versus
topological groups: dichotomy and duality. \textsl{Dissert. Math.} \textbf{%
472} (2010), 138 pp. \newline
\noindent \lbrack BinO3] N. H. Bingham and A. J. Ostaszewski, Regular
variation without limits, \textsl{J. Math. Anal. Appl.} \textbf{370} (2010),
322-338.\newline
\noindent \lbrack BinO4] N. H. Bingham and A. J. Ostaszewski, Dichotomy and
infinite combinatorics: the theorems of Steinhaus and Ostrowski. \textsl{%
Math. Proc. Cambridge Philos. Soc. }\textbf{150} (2011), 1--22.\newline
\noindent \lbrack BinO5] N. H. Bingham and A. J. Ostaszewski, Beurling slow
and regular variation, \textsl{Trans. London Math. Soc., }\textbf{1} (2014)
29-56.\newline
\noindent \lbrack BinO6] N. H. Bingham and A. J. Ostaszewski, Cauchy's
functional equation and extensions: Goldie's equation and inequality, the Go%
\l \k{a}b-Schinzel equation and Beurling's equation. \textsl{Aequationes
Math.} \textbf{89} (2015), 1293--1310.\newline
\noindent \lbrack BinO7] N. H. Bingham and A. J. Ostaszewski, Beurling
moving averages and approximate homomorphisms, \textsl{Indag. Math. }\textbf{%
27} (2016), 601-633 (fuller version: arXiv1407.4093).\newline
\noindent \lbrack BinO8] N. H. Bingham and A. J. Ostaszewski, A. J.
Category-measure duality: convexity, midpoint convexity and Berz
sublinearity. \textsl{Aequationes Math.} \textbf{91} (2017), 801--836.%
\newline
\noindent \lbrack BinO9] N. H. Bingham and A. J. Ostaszewski, Beyond
Lebesgue and Baire IV: Density topologies and a converse Steinhaus-Weil
theorem, \textsl{Topology Appl.}, {\textbf{239} (2018), 274-292 (}%
arXiv1607.00031).\newline
\noindent \lbrack BinO10]{\ N. H. Bingham and A. J. Ostaszewski, Additivity,
subadditivity and linearity: Automatic continuity and quantifier weakening.
\textsl{Indag. Math.} (N.S.) \textbf{29} (2018), 687--713 (arXiv
1405.3948v3).}\newline
\noindent \lbrack BinO11] {N. H. Bingham and A. J. Ostaszewski, Variants on
the Berz sublinearity theorem, }\textsl{Aequationes Math., }Online First,
2018, doi.org/10.1007/s00010-018-0618-8 [arXiv:1712.05183].\newline
\noindent \lbrack BinO12] {N. H. Bingham and A. J. Ostaszewski, Set theory
and the analyst, }\textsl{European. J. Math., }Online First, 2018,
doi.org/10.1007/s40879-018-0278-1 [arXiv:1801.09149].\newline
\noindent \lbrack BinO13] N. H. Bingham and A. J. Ostaszewski, The
Steinhaus-Weil property: its converse, Solecki amenability and
subcontinuity, arXiv:1607.00049v3.\newline
\noindent \lbrack BinO14] {N. H. Bingham and A. J. Ostaszewski, }Sequential
regular variation: extensions of Kendall's theorem, preprint.\newline
\noindent \lbrack Bir] G. Birkhoff, A note on topological groups. \textsl{%
Compos. Math.} \textbf{3} (1936), 427--430.\newline
\noindent \lbrack BojK] R. Bojani\'{c} and J. Karamata, \textsl{On a class
of functions of regular asymptotic behavior, }Math. Research Center Tech.
Report 436, Madison, Wis. 1963; reprinted in \textsl{Selected papers of
Jovan Karamata} (ed. V. Mari\'{c}, Zevod za Ud\v{z}benike, Beograd, 2009),
545-569.\newline
\noindent \lbrack Brz1] J. Brzd\k{e}k, The Go\l \k{a}b-Schinzel equation and
its generalizations, \textsl{Aequat. Math.} \textbf{70} (2005), 14-24.%
\newline
\noindent \lbrack Brz2] J. Brzd\k{e}k, Subgroups of the group $\mathbb{Z}%
_{n} $ and a generalization of the Go\l \k{a}b-Schinzel functional equation,
\textsl{Aequat. Math.} \textbf{43} (1992), 59-71.\newline
\noindent \lbrack BrzM] J. Brzd\k{e}k and A. Mure\'{n}ko, On a conditional Go%
\l \k{a}b-Schinzel equation, \textsl{Arch. Math.} \textbf{84} (2005),
503-511.\newline
\noindent \lbrack Chu1] J. Chudziak, Semigroup-valued solutions of the Go\l
\k{a}b-Schinzel type functional equation, \textsl{Abh. Math. Sem. Univ.
Hamburg,} \textbf{76} (2006), 91-98.

\noindent \lbrack Chu2] J. Chudziak, Semigroup-valued solutions of some
composite equations, \textsl{Aequat. Math.} \textbf{88} (2014), 183-198.%
\newline
\noindent \lbrack CrnGH] M. Crnjac, B. Gulja\v{s} and H. I. Miller, On some
questions of Ger, Grubb and Kraljevi\'{c}. \textsl{Acta Math. Hungar.}
\textbf{57} (1991), 253--257.\newline
\noindent \lbrack Dar] U. B. Darji, On Haar meager sets. \textsl{Topology
Appl.} \textbf{160} (2013), 2396--2400.\newline
\noindent \lbrack DieS] J. Diestel and A. Spalsbury,\textsl{\ The joys of
Haar measure.} Graduate Studies in Mathematics \textbf{150}, Amer. Math.
Soc., 2014.\newline
\noindent \lbrack Ell] R. Ellis, \textsl{Lectures on topological dynamics.}
W. A. Benjamin, 1969.\newline
\noindent \lbrack EllE] D. B. Ellis and R. Ellis, \textsl{Automorphisms and
equivalence relations in topological dynamics.} London Math. Soc. Lect. Note
Series \textbf{412}. Cambridge University Press, 2014.\newline
\noindent \lbrack HarW] G. H. Hardy and E. M. Wright, \textsl{An
introduction to the theory of numbers}. 6$^{\text{th}}$ ed. Revised by D. R.
Heath-Brown and J. H. Silverman. With a foreword by Andrew Wiles. Oxford
University Press, 2008.\newline
\noindent \lbrack Jab1] E. Jab\l o\'{n}ska, Continuous on rays solutions of
an equation of the Go\l \c{a}b-Schinzel type. \textsl{J. Math. Anal. Appl.}
\textbf{375} (2011), 223--229.\newline
\noindent \lbrack Jab2] E. Jab\l o\'{n}ska, On continuous solutions of an
equation of the Go\l \c{a}b-Schinzel type. \textsl{Bull. Aust. Math. Soc.}%
\textbf{\ 87} (2013), 10--17.\newline
\noindent \lbrack Jab3] E. Jab\l o\'{n}ska, On locally bounded above
solutions of an equation of the Go\l \c{a}b-Schinzel type. \textsl{Aequat.
Math.} \textbf{87} (2014),125--133.\newline
\noindent \lbrack Jab4] E. Jab\l o\'{n}ska, On continuous on rays solutions
of a composite-type equation. \textsl{Aequat. Math.} \textbf{89} (2015),
583--590.\newline
\noindent \lbrack Jab5] E. Jab\l o\'{n}ska, On solutions of some
generalizations of the Go\l \c{a}b-Schinzel equation, in: \textsl{Functional
equations in mathematical analysis}, ed. J. Brzd\k{e}k and Th. M. Rassias,
Springer Optim. Appl. \textbf{52} (2012), 509-521.\newline
\noindent \lbrack Jav] P. Javor, On the general solution of the functional
equation $f(x+yf(x))=f(x)f(y).$ \textsl{Aequat. Math.} \textbf{1} (1968),
235-238.\newline
\noindent \lbrack Kak1] S. Kakutani, \"{U}ber die Metrisation der
topologischen Gruppen, \textsl{Proc. Imp. Acad. Tokyo} \textbf{12} (1936),
82--84 (reprinted in [Kak2]).\newline
\noindent \lbrack Kak2] S. Kakutani, \textsl{Selected papers.} Vol. 1. (Ed.
R. R. Kallman), Contemporary Mathematicians. Birkh\"{a}user, 1986.\newline
\noindent \lbrack Kuc] Marek Kuczma, \textsl{An introduction to the theory
of functional equations and inequalities. Cauchy's equation and Jensen's
inequality.} 2nd ed., Birkh\"{a}user, 2009 [1st ed. PWN, Warszawa, 1985].%
\newline
\noindent \lbrack Loo] L. H. Loomis, \textsl{An introduction to abstract
harmonic analysis.} Van Nostrand, 1953.\newline
\noindent \lbrack LukMZ] J. Luke\v{s}, J. Mal\'{y}, L. Zaj\'{\i}\v{c}ek,
\textsl{Fine topology methods in real analysis and potential theory},
Lecture Notes in Mathematics \textbf{1189}, Springer, 1986.\newline
\noindent \lbrack MatZ] E. Mato\u{u}skov\'{a}, M. Zelen\'{y}, A note on
intersections of non--Haar null sets, \textsl{Colloq. Math.} \textbf{96}
(2003), 1-4.\newline
\noindent \lbrack MilMO] H. I. Miller, L. Miller-van Wieren, and A. J.
Ostaszewki, Beyond Erd\H{o}s-Kunen-Mauldin: Singular sets with
shift-compactness properties, preprint.\newline
\noindent \lbrack MilO] H. I. Miller and A.J. Ostaszewski, Group action and
shift-compactness, \textsl{J. Math. Anal. App.} \textbf{392} (2012), 23--39.%
\newline
\noindent \lbrack Ost1] A. J. Ostaszewski, Beurling regular variation, Bloom
dichotomy, and the Go\l \k{a}b-Schinzel functional equation, \textsl{Aequat.
Math.} \textbf{89} (2015), 725-744. \newline
\noindent \lbrack Ost2] A. J. Ostaszewski, Homomorphisms from Functional
Equations: The Goldie Equation, \textsl{Aequat. Math. }\textbf{90} (2016),
427-448 (arXiv: 1407.4089).\newline
\noindent \lbrack Ost3] A. J. Ostaszewski, Homomorphisms from Functional
Equations in Probability, in: \textsl{Developments in Functional Equations
and Related Topics}, ed. J. Brzd\k{e}k et al, Springer Optim. Appl. 124
(2017), 171-213.\newline
\noindent \lbrack Oxt] J. C. Oxtoby: \textsl{Measure and category}, 2nd ed.
Graduate Texts in Math. \textbf{2}, Springer, 1980.\newline
\noindent \lbrack Pet] {B. J. Pettis, {On continuity and openness of
homomorphisms in topological groups.} \textsl{Ann. of Math. }(2) \textbf{52}
(1950), 293--308.}\newline
\noindent \lbrack Pic] {S. Piccard, {Sur les ensembles de distances des
ensembles de points d'un espace Euclidien.\ }\textsl{M\'{e}m. Univ. Neuch%
\^{a}tel} \textbf{13}, 212 pp. 1939.}\newline
\noindent \lbrack Pop] C. G. Popa, Sur l'\'{e}quation fonctionelle $%
f[x+yf(x)]=f(x)f(y),$ \textsl{Ann. Polon. Math.} \textbf{17} (1965), 193-198.%
\newline
\noindent \lbrack Rud1] W. Rudin, \textsl{Fourier analysis on groups}.
Wiley, 1962.\newline
\noindent \lbrack Rud2] W. Rudin, \textsl{Real and complex analysis}. 3$^{%
\text{rd}}$ ed. McGraw-Hill, 1987.\newline
\noindent \lbrack Sak] S. Saks, \textsl{Theory of the integral}, Dover, 1964
(transl. of Mono. Mat. VII, 1937).\newline
\noindent \lbrack Sol] {\ S. Solecki, {Amenability, free subgroups, and Haar
null sets in non-locally compact groups}. \textsl{Proc. London Math. Soc.}
(3) \textbf{93} (2006), 693--722.}\newline
\noindent \lbrack Ste] H. Steinhaus, Sur les distances des points de mesure
positive.\ \textsl{Fund. Math. }\textbf{1} (1920), 83-104.\newline
\noindent \lbrack Wei] A. Weil,\ \textsl{L'int\'{e}gration dans les groupes
topologiques}, Actualit\'{e}s Scientifiques et Industrielles 1145, Hermann,
1965 (1$^{\text{st }}$\ ed. 1940).\newline

\bigskip

\noindent Mathematics Department, Imperial College, London SW7 2AZ;
n.bingham@ic.ac.uk \newline
Mathematics Department, London School of Economics, Houghton Street, London
WC2A 2AE; A.J.Ostaszewski@lse.ac.uk\newpage

\end{document}